\documentclass[12pt,reqno]{amsart}
\usepackage{amsfonts,enumitem}
\usepackage{ifthen}
\usepackage{amsthm}
\usepackage{amsmath}
\usepackage{graphicx}
\usepackage{caption}
\usepackage{subcaption}
\usepackage{amscd,amssymb}
\usepackage{color}
\usepackage{hyperref}
\usepackage{array,multirow,multicol,makecell}

\usepackage{cite,epstopdf}

\setcellgapes{1pt}
\makegapedcells
\newcolumntype{R}[1]{>{\raggedleft\arraybackslash }b{#1}}
\newcolumntype{L}[1]{>{\raggedright\arraybackslash }b{#1}}
\newcolumntype{C}[1]{>{\centering\arraybackslash }b{#1}}
\newcounter{minutes}
\divide\time by 60
\newcounter{hours}
\multiply\time by 60 \addtocounter{minutes}{-\time}

\newtheorem{theorem}{Theorem}
\newtheorem{lemma}{Lemma}
\newtheorem{definition}{Definition}
\newtheorem{corollary}{Corollary}

\newtheorem{example}{Example}

\numberwithin{equation}{section}

\title[Some geometric properties of certain class of Le Roy type functions]{Some geometric properties of certain class of Le Roy type functions}

\author[A. Sharma, S. B. Mahesh, A. Kumar,  K.V. Pai]
{Abhinav Sharma, Suhas B Mahesh, Anish Kumar and Karthik V Pai}

\address{{\bf Abhinav Sharma \footnotemark[1]}\newline
Indian Institute of Technology Hyderabad,\newline
Hyderabad, Telangana, India}
\email{abhisha8055@gmail.com}
\address{{\bf Suhas B Mahesh}\newline
Shiv Nadar University,\newline
Delhi-NCR, India}
\email{suhasbmahesh11@gmail.com}
\address{{\bf Anish Kumar}\newline
Department of Mathematics,
Dr. Shyama Prasad Mukherjee University,\newline
Ranchi 834008, Jharkhand, India}
\email{2019rsma003@nitjsr.ac.in}
\address{{\bf Karthik V Pai \footnotemark[1]\footnotetext[1]{This work was carried out by Abhinav Sharma and Karthik V Pai when they was doing their Masters in Mathematics at RIE Mysore. }}\newline
Indian Institute of Technology Madras,\newline
Chennai, Tamil Nadu, India}
\email{karthikvpai777@gmail.com}

\begin{document}
\maketitle 

\begin{abstract} The main goal of this paper is to obtain sufficient conditions so that Le Roy type functions and multivariate Le Roy type functions satisfy subordination of exponential function.  Moreover conditions on parameters have been derived to claim them being exponential starlike and exponential convex for both of the functions. Starlikeness, convexity and close-to-convexity have also been studied for multivariate Le Roy type functions. Results developed in this work are presumably new and their significance is illustrated by several consequences, graphical representations and examples. 
\end{abstract}

\section{General Introduction and Motivation}
Let $ \mathcal{A}$ denote the class of analytic functions inside the unit disc $\mathcal{D},$ 
$$ \mathcal{D}=\{z\in\mathbb{C}, |z|<1\},$$ having the form \begin{equation}
    f(z)= z +\sum_{k=2}^\infty a_{k}z^{k} ,
\end{equation}

where $a_{k}\in\mathbb{C}$, such that $f(0)$=$f'(0)-1=0.$ This is infinite series representation of analytic functions, where each $a_{k}$ denotes the coefficient of the Taylor's series. The study of geometric properties is reliant on the coefficient conditions and its bounds. This problem is closely related to the famous Bieberbach Conjecture\cite{Bieberbach conjecture}.

Let $S$ denote the set of all functions that are univalent about the origin in unit disk $\mathcal{D}$. A univalent function is said to be starlike in $\mathcal{D}$ if $f \in \mathcal{A}$ and $f$ satisfies the condition  $Re\left({\frac{zf^{\prime}(z)}{f(z)}}\right)>0$, which is denoted as $S\mbox{*}$ . That is
$$f \in S\mbox{*} \iff Re\left({\frac{zf^{\prime}(z)}{f(z)}}\right)>0.$$
Let $0 \le\ \alpha\ < 1$. Then a function $f$ is said to be starlike of order $\alpha$ if $Re\left({\frac{zf^{\prime}(z)}{f(z)}}\right)>\alpha$.
i.e  $$f \in S\mbox{*}(\alpha) \iff Re\left({\frac{zf^{\prime}(z)}{f(z)}}\right)>\alpha.$$

Similarly, the class of convexity of univalent functions denoted by $K\mbox{*}$ is defined by its analytic characterization as follows:

  $$f \in K\mbox{*} \iff Re\left({1+\frac{zf^{\prime\prime}(z)}{f^{\prime}(z)}}\right)>0, $$

where   $f \in \mathcal{A}\ $.
\newpage
A function  $f \in \mathcal{A}\ $ is said to be convex of order $\alpha$ if $0\le \alpha$ $<1$
and $f$ satisfies the following criteria

  $$f \in K\mbox{*}(\alpha) \iff Re\left({1+\frac{zf^{\prime\prime}(z)}{f^{\prime}(z)}}\right)>\alpha.$$

An analytic function $g$ in $\mathcal{A}$ is called close-to-convex
in the open unit disk $\mathcal{D}$ if there exists a function
$h(z)$, which is starlike in $\mathcal{D}$ such that
$$Re\left(\dfrac{zg'(z)}{h(z)}\right)>0\qquad (z\in \mathcal{D}).$$

The concept of subordination plays a crucial role in the classification of various classes of analytic functions on the open unit disk $\mathcal{D}$. An analytic function $g(z)$ is said to be subordinate to analytic and univalent function $f(z)$ if $g(0)=f(0)$ for all $z$ in $\mathcal{D}$ and the image set of $g(z)$ is contained in image set of $f(z)$.

      That is $$g(z) \subset f(z),  \forall z \in \mathcal{D}.$$

Let us consider ${P_e}$ be the subclass of all univalent functions, which is defined as the class of all analytic functions satisfy the  criteria $f(0)=1$ and $f(z) \subset e^{z}$ on the $\mathcal{D}.$ Now we define exponential convexity and starlikenes as, if the function $f \in \mathcal{A}$ and the quantities $\left(1+\frac{zf^{\prime\prime}(z)}{f'(z)}\right)$ and $\left(\frac{zf^{\prime}(z)}{f(z)}\right)$ be such that it belongs to ${P_e}$, then the function $f$ is exponentially convex and exponentially starlike respectively\cite{Geometric Properties}. We denote such functions as ${K_e}$ and ${S_e}$ respectively. These functions are closely related to the subclasses of convex and starlike functions introduced by Ma and Minda\cite{Geometric Properties}.


Mittag-Leffler and Le Roy type functions are two types of special functions that arise frequently in various branches of mathematics and physics. The Mittag-Leffler function \cite{Mittag Leffler} is a generalization of the exponential function, which is known as the queen function of fractional calculus, defined by an infinite series involving the gamma function. It is denoted by $E_{\alpha, \beta}(z)$ and depends on two parameters $\alpha$ and $\beta$. It is defined as

$$  E_{\alpha,\beta}(z) = \sum_{n=0}^{\infty} \frac{z^n}{\Gamma(\alpha n + \beta)}, $$ where $\alpha$ in $\mathcal{C}$ and $\beta \in \mathcal{C}$. 
 This complex variable function turns out be entire under the conditions $\emph{Re}(\alpha)>0$ and $\beta>0.$
 The function has various applications in fractional calculus, anomalous diffusion, and in the study of viscoelastic materials.

On the other hand, the Le Roy type function is a multivariate function that is defined by the sum of power functions raised to their own index. It is represented by $$  F_{1,1}^{\gamma}(z) = \sum_{n=0}^{\infty} \frac{z^n}{{(n!)}^{\gamma}}. $$ The function is useful in the study of classical mechanics, fluid dynamics, and turbulence.

Both of these functions have important mathematical properties and applications in diverse fields of science and engineering. They have been extensively studied by researchers and have led to many interesting results and discoveries.
    We now define the special function Le Roy type Mittag Leffler Function
  $$ F_{1,1}^{\gamma}(z) = \sum_{n=0}^{\infty} \frac{z^n}{(n!)^\gamma} .$$

  We can clearly see that the Le Roy type function is a generalization of Mittag Leffler function or exponential function which can also be expressed as
  $$F_{1,1}^{\gamma}(z) = \sum_{n=0}^{\infty} \frac{z^n}{[\Gamma(n+1)]^\gamma},$$
  
  $$  F_{\alpha,\beta}^{\gamma}(z) = \sum_{n=0}^{\infty} \frac{z^n}{[\Gamma(\alpha n + \beta)]^\gamma}. $$
    It turns out that this complex function of three parameter is entire over the complex plane under the assumptions that $\mathbb{Re}(\alpha)>0$ and $\beta > 0$ and $\gamma>0$ and R Garrappa in his paper\cite{Le Roy type} on the integral transforms of the 3 parameter Wright functions, studied the growth and order of the Le Roy functions and provided the integral transforms.
  Le Roy and Mittag Leffler were French Mathematicians, who were working in competition to develop functions as solutions for fractional differential equations. Le Roy type Mittag Leffler Functions are the generalizations of Mittag Leffler functions of two parameters\cite{Le Roy type}:
  $$  F_{\alpha,\beta}^{\gamma}(z) = \sum_{n=0}^{\infty} \frac{z^n}{[\Gamma(\alpha n + \beta)]^\gamma},$$
  where $\Gamma(\alpha n + \beta)$ is the generalized gamma function over the Complex plane defined by \cite{Textbook}
  $$ \Gamma(z) = \begin{cases}
    \frac{\pi}{\sin(\pi z)\Gamma(1-z)} & \text{if }\Re(z) \le 0 \\
    \Gamma(z) & \text{otherwise},
\end{cases}
$$
where 
$$\Gamma (z) = \int_{0}^{\infty} t^{z-1} e^{-t}dt.$$
Upon substituting specific values for each of the parameters  
$(\alpha, \beta, \gamma)$ we have various prominent special functions, for example
\begin{enumerate}
    \item 
  $F_{1,1}^{1}(z) = \sum_{n=0}^{\infty} \frac{z^n}{[\Gamma(n + 1)]} = e^z$
\item 
$F_{2,2}^{1}(z) = \frac{\sinh{\sqrt{z}}}{z}$

\item 
  $  F_{\alpha,\beta}^{1}(z) = \sum_{n=0}^{\infty} \frac{z^n}{[\Gamma(\alpha n + \beta)]} = E_{\alpha, \beta}(z)$
\item 
    $ F_{1,1}^{2}(z) = J_0(\sqrt{z}) $
\item 
$  F_{\alpha,\beta}^{n}(z) = \sum_{n=0}^{\infty} \frac{z^n}{[\Gamma(\alpha n + \beta)]^n} = E^{n}_{(\alpha, \beta),...(\alpha,\beta)}(z)$
\item 
$ R^{\gamma}_{1,1} = \sum_{n=0}^{\infty} \frac{z^n}{{(k!)}^{\gamma}} = F_{1,1}^{\gamma}(z) $
\end{enumerate}

Now we extend the definition of the generalization of Le Roy functions of three parameters to 3n parameters over the complex domain. The extension is quite natural using the gamma functions and finite number of parameters, we identify the indices with natural indexes\cite{Multiindex Le Roy}:

The multi-index function is defined by
$$F^{\gamma_{1},\gamma_{2},\ldots,\gamma_{n}}_{\alpha_{1},\alpha_{2},\ldots,\alpha_{n},\beta_{1},\beta_{2},\ldots,\beta_{n}}(z)
= \sum_{m=0}^{\infty} \frac{z^m}{\prod_{i=1}^{n}[\Gamma(\alpha_i m + \beta_i)]^{\gamma_i}},$$
and, for brevity, we write
$$\mathcal{F}(z):=F^{\gamma_{1},\gamma_{2},\ldots,\gamma_{n}}_{\alpha_{1},\alpha_{2},\ldots,\alpha_{n},\beta_{1},\beta_{2},\ldots,\beta_{n}}(z).$$
Its standard normalized form is defined by
$$\mathbb{F}^{\gamma_{1},\gamma_{2},\ldots,\gamma_{n}}_{\alpha_{1},\alpha_{2},\ldots,\alpha_{n},\beta_{1},\beta_{2},\ldots,\beta_{n}}(z)
:=z\prod_{i=1}^{n}[\Gamma(\beta_i)]^{\gamma_i}\mathcal{F}(z),$$
and we denote it by
$$\mathbf{F}(z):=\mathbb{F}^{\gamma_{1},\gamma_{2},\ldots,\gamma_{n}}_{\alpha_{1},\alpha_{2},\ldots,\alpha_{n},\beta_{1},\beta_{2},\ldots,\beta_{n}}(z).$$
The Le Roy type Mittag-Leffler function is clearly a special case of the above mentioned function and it can be obtained in many combinations of the parameters.

The study of geometric properties of special functions has been considered a vital area of research for a long period of time, several mathematicians have worked in this regard. There has been various works carried out on geometric properties of special functions such as Bessel functions, Hypergeometric function, Fox-Wright functions, Mittag-Leffler functions etc \cite{kumar,  Medirata, MehrezBansal, bansla}. With these insightful works as the background we proceed in this investigation to study few important geometric properties of the Le Roy type Mittag-Leffler and Multivariate Le Roy type functions. In the next section we consider some useful definitions and lemmas to derive our main results.
\section{Useful definitions and lemmas}



\begin{definition} \cite{Subordination Definition}
     An analytic function $f(z)$ is said to be the subordinate of analytic and univalent function $g(z)$ on the unit disk $|z| < 1$, if there exists a schwarz function $w(z)$ such that $w(0)= 0$ with $|w(z)| \le |z| <1 $ and $f(z)=g(w(z))$ in $|z| < 1.$
\end{definition}


\begin{lemma}\cite{Convexity and Starlikeness}\label{Starlike}
Let $f(z) \in \mathcal{A}$ and $|f'(z)-1| < \frac{2}{\sqrt{5}}$ $\forall z \in \mathcal{D}, $ then f(z) is starlike in $\mathcal{D}.$
\end{lemma} \begin{lemma}\cite{Medirata}\label{Subordination Lemma}
Let $S_{1}$ be the image set of the function $f(z)$ and $S_{2}$ be the image set of the function $e^z$ then $S_{1} \subset  S_{2}$ if $$|f(z)-1| < 1-\frac{1}{e},$$ and $f(0)=1$,
 where e is the Euler's number.
\end{lemma}
\begin{lemma}\cite{mehrez}\label{Digamma Inequality}
    For any positive real number s $\ge$ 1, the digamma function (psi-function) $\psi(s)=\frac{\Gamma'(s)}{\Gamma(s)}$ satisfies the following inequality 
    
    $$\log(s)- \delta \leq \psi(s)< \log(s),$$
    where $\delta$ is the Euler Mascheroni constant.
\end{lemma}
\begin{lemma}[Ozaki's Lemma]\cite{Ozaki}
    Let $f(z)=z+\sum_{k=2}^{\infty} A_{k} z^{k}$. If $1 \leq 2 A_{2} \leq \ldots \leq n A_{n} \leq$ $(n+1) A_{n+1} \leq \ldots \leq 2$, or $1 \geq 2 A_{2} \geq \ldots \geq n A_{n} \geq(n+1) A_{n+1} \geq \ldots \geq 0$, then $f$ is close-to-convex with respect to $-\log (1-z)$.
\end{lemma}

\begin{lemma}\cite{mcgregor}\label{mcgregorlemma}
    Let $f \in \mathcal{A}$ and $|(f(z) / z)-1|<1$ for each $z \in \mathcal{D}$, then $f$ is univalent and starlike in $\mathcal{D}_{1 / 2}=\{z:|z|<1 / 2\}$.
\end{lemma}
\begin{lemma}\label{convexnew}\cite{thmac}
    Let $f \in \mathcal{A}$ and $\left|f^{\prime}(z)-1\right|<1$ for each $z \in \mathcal{D}$, then $f$ is convex in $\mathcal{D}_{1 / 2}=\{z:|z|<1 / 2\}$.
\end{lemma}

\begin{lemma}\label{newlem}
     Let
     $$\mathbf{F}(z)=\mathbb{F}^{\gamma_{1},\gamma_{2},\ldots,\gamma_{n}}_{\alpha_{1},\alpha_{2},\ldots,\alpha_{n},\beta_{1},\beta_{2},\ldots,\beta_{n}}(z)$$
     be the standard normalized multi-index Le Roy type function. Let $\min (\alpha_i, \gamma_i) \geq 1, \beta_i>0$ for all $i=1,2,\cdots,n$ such that $\alpha_j+\beta_j \geq 2$ for some $j \in {1,2,\cdots,n}$. Then the following inequality
\begin{equation*}
\mathbf{F}(z) \leq z+z \theta^{\gamma_{1},\gamma_{2},\ldots,\gamma_{n}}_{\alpha_{1},\alpha_{2},\ldots,\alpha_{n},\beta_{1},\beta_{2},\ldots, \beta_{n}}\left(e^{z}-1\right) 
\end{equation*}
holds true for all $z>0$, where

\begin{equation*}
\theta^{\gamma_{1},\gamma_{2},\ldots,\gamma_{n}}_{\alpha_{1},\alpha_{2},\ldots,\alpha_{n},\beta_{1},\beta_{2},\ldots, \beta_{n}}=\prod_{i=1}^{n}\left[\frac{\Gamma(\beta_i)}{\Gamma(\alpha_{i}+\beta_{i})}\right]^{\gamma_i}
\end{equation*}
\end{lemma}
\begin{proof}
We need to prove that the sequence
\begin{equation*}
x_{k}:=\left\{\frac{\Gamma(k+1)}{\prod_{i=1}^{n}[\Gamma(\alpha_{i}k+\beta_{i})]^{\gamma_{i}}}\right\}_{k \geq 1} 
\end{equation*}
is decreasing. Let $\min (\alpha_i, \gamma_i) \geq 1$ and $\beta_i>0$ for all $i=1,2,\cdots,n$, then we have

\begin{align*}
\frac{x_{k+1}}{x_{k}} & =\frac{(k+1)\prod_{i=1}^{n}[\Gamma(\alpha_i k+\beta_i)]^{\gamma_i}}{\prod_{i=1}^{n}[\Gamma(\alpha_i k+\alpha_i+\beta_i)]^{\gamma_i}} \\
& \leq \frac{(k+1)\prod_{i=1}^{n}[\Gamma(\alpha_i k+\beta_i)]^{\gamma_i}}{\prod_{i=1}^{n}[\Gamma(\alpha_i k+\beta_i+1)]^{\gamma_i}}=\frac{k+1}{\prod_{i=1}^{n}(\alpha_i k+\beta_i)^{\gamma_i}} \leq \frac{k+1}{(\alpha_j k+\beta_j)} 
\end{align*}
for some $j \in {1,2,\cdots,n}$. The first inequality uses the fact that $\Gamma$ is increasing on $[2,\infty)$ and $\alpha_i k+\beta_i+1\ge\alpha_i+\beta_i+1\ge 2$, so $\Gamma(\alpha_i k+\alpha_i+\beta_i)\ge\Gamma(\alpha_i k+\beta_i+1)$ since $\alpha_i\ge 1$. The last inequality follows because $\alpha_i k+\beta_i\ge 1$ and $\gamma_i\ge 1$ imply $(\alpha_i k+\beta_i)^{\gamma_i}\ge 1$, hence the product exceeds any single factor $\alpha_j k+\beta_j$.
Thus $\frac{x_{k+1}}{x_{k}} \le \frac{k+1}{\alpha_j k+\beta_j} \le 1$ exactly when $k+1\le \alpha_j k+\beta_j$, i.e.\ when
$$\chi(k):=(\alpha_j-1)k+\beta_j-1 = (\alpha_j-1)(k-1)+(\alpha_j+\beta_j)-2\ge 0.$$
This in turn implies that the sequence $\left(x_{k}\right)_{k \geq 1}$ monotonically decreases. Therefore, for $z>0$ we get

$$
\begin{aligned}
\frac{\mathbf{F}(z)}{z} & =1+\sum_{k=1}^{\infty} \frac{\prod_{i=1}^{n}[\Gamma(\beta_i)]^{\gamma_i} \Gamma(k+1)}{\prod_{i=1}^{n}[\Gamma(\alpha_i k+\beta_i)]^{\gamma_i}} \frac{z^{k}}{k!} \\
& \leq 1+\theta^{\gamma_{1},\gamma_{2},\ldots,\gamma_{n}}_{\alpha_{1},\alpha_{2},\ldots,\alpha_{n},\beta_{1},\beta_{2},\ldots, \beta_{n}} \sum_{k=1}^{\infty} \frac{z^{k}}{k!}=1+\theta^{\gamma_{1},\gamma_{2},\ldots,\gamma_{n}}_{\alpha_{1},\alpha_{2},\ldots,\alpha_{n},\beta_{1},\beta_{2},\ldots, \beta_{n}}\left(e^{z}-1\right).
\end{aligned}
$$
\end{proof}
\section{exponential subordination}
 In this section, we consider the family of Le Roy type Mittag- Leffler functions, multi index Le Roy type functions and study some of its geometric properties of the family of functions which satisfy the criteria of exponential subordination. 
\begin{theorem}\label{Le Roy Subordination}
  The Le Roy type Mittag-Leffler function $F_{\alpha,\beta}^{\gamma}$ is said to be subordinate of exponential function $e^{z}$ on the $\mathcal{D}$ under the condition that $\alpha, \beta, \gamma>0$, $\alpha\gamma \ge 1$, $\alpha^{2}\gamma \ge \beta$, ${[\Gamma(\beta)]^\gamma}=1$, $\log(2)-\alpha\gamma \log(\alpha  + \beta) + \alpha\gamma\delta<0$ and $$[\Gamma(\alpha + \beta)]^{\gamma} > e.$$
\end{theorem}
\begin{proof}
To prove the above theorem, we will use the Lemma \ref{Subordination Lemma} which gives us a sufficient condition to guarantee that a function $f(z)$ is contained in the image set of exponential function.
Consider the Le Roy type function $F_{\alpha,\beta}^{\gamma}(z)=\sum_{n=0}^{\infty} \frac{z^n}{[\Gamma(\alpha n + \beta)]^\gamma},$ which satisfy

$$F_{\alpha,\beta}^{\gamma}(0)=\frac{1}{[\Gamma(\beta)]^\gamma}=1.$$
Now, it suffices to show that $$|F_{\alpha,\beta}^{\gamma}(z)-1| < 1-\frac{1}{e}.$$ 
We have,
$$ \left|\sum_{n=0}^{\infty} \frac{z^n}{[\Gamma(\alpha n + \beta)]^\gamma}-1\right| = \left|\sum_{n=0}^{\infty} \frac{\Gamma(n+1)z^n}{n![\Gamma(\alpha n + \beta)]^\gamma}-1\right|. $$
Define the coefficients 
 $$b_{n}(\alpha,\beta,\gamma)=\frac{\Gamma(n+1)}{[\Gamma(\alpha n + \beta)]^\gamma},\qquad n\ge 1.$$
To show that $\{b_n\}_{n\ge 1}$ is decreasing, consider the function $f_{1}(s)$ defined by  $$f_{1}(s)=\frac{\Gamma(s+1)}{[\Gamma(\alpha s+\beta)]^{\gamma}},  s\ge1,$$
$$ \implies f'_{1}(s)=f_{1}(s)f_{2}(s),$$
 where $$f_{2}(s)=\psi(s+1)-\alpha\gamma\psi(\alpha s+\beta),\qquad s\ge 1.$$
The hypothesis
\[
\log2-\alpha\gamma\log(\alpha+\beta)+\alpha\gamma\delta<0
\implies
\alpha\gamma\log(\alpha+\beta)
>\log2+\alpha\gamma\delta>0.
\]
Therefore $\alpha+\beta>1$, and for every $s\ge1$, $
\alpha s+\beta\ge\alpha+\beta>1.
$
So the digamma bounds in Lemma~\ref{Digamma Inequality} can be applied to both $s+1$ and $\alpha s+\beta$. Thus, using Lemma \ref{Digamma Inequality} to bound the digamma functions, we define $f_{3}(s)$ such that $f_{2}(s) < f_{3}(s)$, namely
$$f_{3}(s) := \log(s+1)-\alpha\gamma \log(\alpha s + \beta) + \alpha\gamma\delta.$$
Upon evaluating and simplification of derivative we get the expression
$$f'_{3}(s) = \frac{\alpha(1-\alpha\gamma)s+(\beta-\alpha^2\gamma)}{(\alpha s + \beta)(s+1)}$$\label{Equation 1}
 By the derivative test, $f_{3}(s)$ is decreasing on $[1,\infty)$ under the conditions $\alpha\gamma \ge 1$ and $\alpha^2\gamma \ge \beta$ (since its numerator is $\alpha(1-\alpha\gamma)s+(\beta-\alpha^2\gamma)\le 0$). Moreover, $f_{3}(1)<0$ by the hypothesis $\log 2-\alpha\gamma\log(\alpha+\beta)+\alpha\gamma\delta<0$. Consequently $f_{3}(s)\le f_{3}(1)<0$ for all $s\ge 1$, hence $f_{2}(s)<0$ and $f_{1}(s)$ is decreasing on $[1,\infty)$. Thus the coefficients $b_{n}(\alpha,\beta,\gamma)$ form a decreasing sequence. 
A simple computation gives
\begin{align*}
&\left|\sum_{n=0}^{\infty} \frac{\Gamma(n+1)z^n}{n![\Gamma(\alpha n + \beta)]^\gamma}-1\right| \leq \sum_{n=1}^{\infty} \left|\frac{b_{1}(\alpha,\beta,\gamma)}{n!}\right|\\
& = b_{1}(\alpha,\beta,\gamma)(e-1) = \left|\frac{e-1}{[\Gamma(\alpha + \beta)]^\gamma}\right| < 1-\frac{1}{e}.
\end{align*}
$$\iff [\Gamma(\alpha + \beta)]^{\gamma} > e,$$
which is the desired condition and concludes the proof.
\end{proof}

\begin{example}
Let us consider the parameters $\alpha=2.5$, $\beta=1$, and $\gamma=2$. This set of values satisfies all the hypotheses of Theorem \ref{Le Roy Subordination}. In Figure~\ref{thm1}, the red region is the image of the unit disk $\mathcal{D}$ under $F_{\alpha,\beta}^{\gamma}(z)$, and the yellow background is the image of $\mathcal{D}$ under $e^{z}$. The red region lies entirely inside the yellow, confirming that $F_{\alpha,\beta}^{\gamma}(z)\subset e^{z}$ on $\mathcal{D}$.
\begin{figure}[h!]
    \centering
    \includegraphics[width=0.5\linewidth]{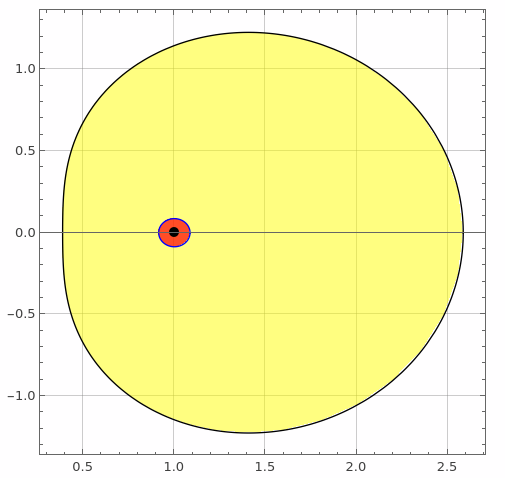}
    \caption{Subordination of $F_{2.5,1}^{2}(z)$ (red) to $e^{z}$ (yellow) on $\mathcal{D}$}
    \label{thm1}
\end{figure}

\end{example}
In the next theorem we will try to obtain sufficient conditions for the exponential subordination of the multi index Le Roy type functions.

\begin{theorem}\label{thm:generalized-leroy}
    Let $\mathcal{F}(z)=F^{\gamma_{1},\gamma_{2},\ldots,\gamma_{n}}_{\alpha_{1},\alpha_{2},\ldots,\alpha_{n},\beta_{1},\beta_{2},\ldots,\beta_{n}}(z) = \sum_{m=0}^{\infty}\frac{z^m}{\prod_{i=1}^{n}[\Gamma(\alpha_{i}m+\beta_{i})]^{\gamma_{i}}}$ be the generalized Le Roy type function of $3n$ parameters. Suppose $\min(\alpha_{i},\gamma_{i})\ge 1$ and $\beta_{i}>0$ for all $i=1,\dots,n$, and $\alpha_{j}+\beta_{j}\ge 2$ for some $j\in\{1,\dots,n\}$. If ${\prod_{i=1}^{n}[\Gamma(\beta_{i})]^{\gamma_{i}}}=1$ and
    $${\prod_{i=1}^{n}[\Gamma(\alpha_{i}+\beta_{i})]^{\gamma_{i}}} > e,$$
    then $\mathcal{F}(z)$ is subordinate to the exponential function $e^{z}$ on $\mathcal{D}$.
\end{theorem}
    \begin{proof}
    To prove this theorem, it suffices to show $\left| \mathcal{F}(z) - 1\right| < 1- \frac{1}{e}$, with the help of Lemma~\ref{Subordination Lemma}. First,  $\mathcal{F}(0)=\frac{1}{\prod_{i=1}^{n}[\Gamma(\beta_{i})]^{\gamma_{i}}}=1$. We have,
    \begin{align*}
        \left| \mathcal{F}(z) - 1\right|
        & = \left|\sum_{k=0}^{\infty}\frac{z^{k}}{\prod_{i=1}^{n}[\Gamma(\alpha_{i}k+\beta_{i})]^{\gamma_{i}}} - 1 \right|
        = \left| \sum_{k=1}^{\infty}\frac{z^{k}}{\prod_{i=1}^{n}[\Gamma(\alpha_{i}k+\beta_{i})]^{\gamma_{i}}} \right| \\
        & = \left| \sum_{k=1}^{\infty}\frac{\Gamma(k+1)\,z^{k}}{k!\prod_{i=1}^{n}[\Gamma(\alpha_{i}k+\beta_{i})]^{\gamma_{i}}} \right|
    \end{align*}
    Since $z$ belongs to the open unit disk,
    $$ \left| \sum_{k=1}^{\infty}\frac{\Gamma(k+1)z^{k}}{k!\prod_{i=1}^{n}[\Gamma(\alpha_{i}k+\beta_{i})]^{\gamma_{i}}} \right|
        \le \sum_{k=1}^{\infty}\frac{\Gamma(k+1)}{k!\prod_{i=1}^{n}[\Gamma(\alpha_{i}k+\beta_{i})]^{\gamma_{i}}}. $$
    Set $x_{k}:=\frac{\Gamma(k+1)}{\prod_{i=1}^{n}[\Gamma(\alpha_{i}k+\beta_{i})]^{\gamma_{i}}}$ for $k\ge 1$. By Lemma~\ref{newlem} the sequence $\{x_{k}\}_{k\ge 1}$ is decreasing under the given hypotheses, hence $x_{k}\le x_{1}= \frac{1}{\prod_{i=1}^{n}[\Gamma(\alpha_{i}+\beta_{i})]^{\gamma_{i}}}$ for all $k\ge 1$. Therefore
    \begin{align*}
        \left| \mathcal{F}(z) - 1 \right|
        & \le \sum_{k=1}^{\infty}\frac{x_{k}}{k!}
        \le x_{1}\sum_{k=1}^{\infty}\frac{1}{k!}
        = \frac{e-1}{\prod_{i=1}^{n}[\Gamma(\alpha_{i}+\beta_{i})]^{\gamma_{i}}}.
    \end{align*}
    The condition $\prod_{i=1}^{n}[\Gamma(\alpha_{i}+\beta_{i})]^{\gamma_{i}} > e$ is equivalent to $\frac{e-1}{\prod_{i=1}^{n}[\Gamma(\alpha_{i}+\beta_{i})]^{\gamma_{i}}} < \frac{e-1}{e}$, and we conclude that $|\mathcal{F}(z)-1| < (e-1)/e$ for all $z\in\mathcal{D}$. By Lemma~\ref{Subordination Lemma}, $\mathcal{F}(z)$ is subordinate to $e^{z}$ on $\mathcal{D}$.
\end{proof}

This theorem is consistent with our previous Theorem\ref{Le Roy Subordination} and thus Le Roy Subordination can be obtained as a corollary of the result on Generalized Le Roy type function of $3n$ parameters.
\begin{example}
Let $n=3$ and consider the parameter sets: $(\alpha_1, \beta_1, \gamma_1) = (2, 1, 2)$,
$(\alpha_2, \beta_2, \gamma_2) = (1, 2, 1)$ and $(\alpha_3, \beta_3, \gamma_3) = (1, 1, 1)$.
Here $\prod_{i=1}^{3}[\Gamma(\beta_i)]^{\gamma_i}=\Gamma(1)^2\cdot\Gamma(2)\cdot\Gamma(1)=1$,
$\min(\alpha_i,\gamma_i)\ge1$, $\alpha_1+\beta_1=3\ge 2$, and
$\prod_{i=1}^{3}[\Gamma(\alpha_i+\beta_i)]^{\gamma_i}= \Gamma(3)^2\cdot\Gamma(3)\cdot\Gamma(2)=8>e$.
Thus all the hypotheses of Theorem~\ref{thm:generalized-leroy} are satisfied. In Figure~\ref{thm2}, the red region is the image of the unit disk $\mathcal{D}$ under the generalized Le Roy type function $\mathcal{F}(z)$, and the yellow background is the image of $\mathcal{D}$ under $e^{z}$. The containment of the red region within the yellow illustrates that $\mathcal{F}(z)\subset e^{z}$ on $\mathcal{D}$.
\begin{figure}[h!]
    \centering
    \includegraphics[width=0.5\linewidth]{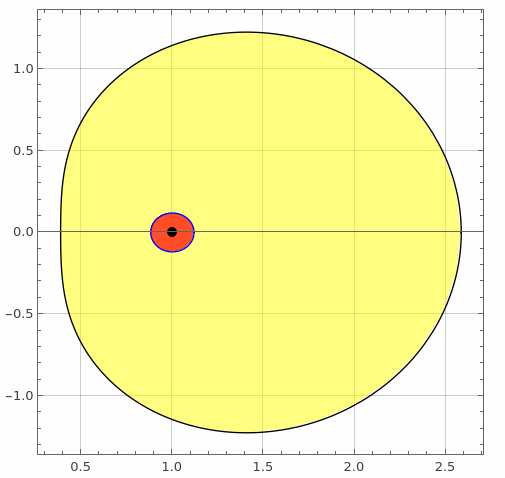}
    \caption{Subordination of the generalized Le Roy type function (red) to $e^{z}$ (yellow) for $n=3$ with $(\alpha_1,\beta_1,\gamma_1)=(2,1,2)$, $(\alpha_2,\beta_2,\gamma_2)=(1,2,1)$,
    $(\alpha_3,\beta_3,\gamma_3)=(1,1,1)$}
    \label{thm2}
\end{figure}
\end{example}

\section{Exponential Convexity and Exponential Starlikeness}
In this section, sufficient conditions for the exponential starlikeness, exponential convexity of Le Roy functions of Mittag Leffler type and multi index Le Roy type functions have been provided using Lemma \ref{Subordination Lemma}. Now, let us consider the Le Roy type Mittag leffler function of three parameters given by the expression
 $$  F_{\alpha,\beta}^{\gamma}(z) = \sum_{n=0}^{\infty} \frac{z^n}{[\Gamma(\alpha n + \beta)]^\gamma}.$$
 Let us consider the following normalization of the above function :
\begin{align*}
&\mathbb{F}_{\alpha,\beta}^{\gamma}(z)= z[\Gamma(\beta)]^{\gamma}F_{\alpha,\beta}^{\gamma}(z) = z+ \sum_{n=2}^{\infty} \frac{z^{n}[\Gamma(\beta)]^{\gamma}}{[\Gamma(\alpha (n-1)+ \beta)]^\gamma} \\
&  {\mathbb{F}_{\alpha,\beta}^{\gamma}}'(z) =1+ \sum_{n=1}^{\infty} \frac{(n+1)z^{n}[\Gamma(\beta)]^{\gamma}}{[\Gamma(\alpha n + \beta)]^\gamma} 
\end{align*}

\begin{theorem}\label{thm:exp-starlike-convex}
      Let $\alpha\ge 1$, $\gamma\ge 1$, $\beta>0$, and $\alpha+\beta\ge 4$.
Then, the normalized Le Roy type function $\mathbb{F}_{\alpha,\beta}^{\gamma}$ is exponentially starlike in $\mathcal{D}$ if it satisfies the condition
  $$ (2e-1)[\Gamma(\beta)]^{\gamma}< [\Gamma(\alpha + \beta)]^{\gamma},$$
 and it is exponentially convex in $\mathcal{D}$ if it satisfies the condition
  $$ (4e^{2}-4e+2)[\Gamma(\beta)]^\gamma < (e-1)[\Gamma(\alpha +\beta)]^\gamma.$$
  \end{theorem}
\begin{proof}
Consider the normalized Le Roy type function
$$\mathbb{F}_{\alpha,\beta}^{\gamma}(z) = z + \sum_{k=2}^{\infty}
\frac{[\Gamma(\beta)]^{\gamma}\,z^{k}}{[\Gamma(\alpha(k-1)+\beta)]^{\gamma}}.$$
Define the base sequence and its weighted counterparts
$$x_{k} = \frac{\Gamma(k+1)}{[\Gamma(\alpha k+\beta)]^{\gamma}},\qquad
w_{k}=k\,x_{k},\qquad z_{k}=(k+1)x_{k},\qquad k\ge 1,$$
and set $C=[\Gamma(\beta)]^{\gamma}$ so that
$c_{n}=C x_{n}$, $b_{n}=C w_{n}$, and $d_{n}=g_{n}=C z_{n}$.\\[1em]
\noindent\textit{Monotonicity.}
Since $\alpha\ge 1$ and $\gamma\ge 1$,
\begin{align*}
\frac{x_{k+1}}{x_{k}}
  &= (k+1)\frac{[\Gamma(\alpha k+\beta)]^{\gamma}}{[\Gamma(\alpha(k+1)+\beta)]^{\gamma}}
   \le \frac{k+1}{(\alpha k+\beta)^{\gamma}}
   \le \frac{k+1}{\alpha k+\beta},
\end{align*}
where we used $\Gamma(x+\alpha)\ge\Gamma(x+1)$ for $x=\alpha k+\beta\ge\alpha+\beta>0$ (this holds because $\Gamma$ is increasing on $[2,\infty)$ and $x+1=\alpha k+\beta+1\ge 2$, while $\alpha\ge1$ ensures $x+\alpha\ge x+1$), together with $(\alpha k+\beta)^{\gamma}\ge\alpha k+\beta$ (as $\alpha k+\beta\ge1$ and $\gamma\ge1$). Hence
\begin{align*}
\frac{w_{k+1}}{w_{k}}
  = \frac{k+1}{k}\,\frac{x_{k+1}}{x_{k}}
  \le \frac{(k+1)^{2}}{k(\alpha k+\beta)},
\qquad
\frac{z_{k+1}}{z_{k}}
  = \frac{k+2}{k+1}\,\frac{x_{k+1}}{x_{k}}
  \le \frac{k+2}{\alpha k+\beta}.
\end{align*}
The upper bound for $w_{k+1}/w_k$ is at most $1$ whenever $$\phi(k):=(\alpha-1)k^2+(\beta-2)k-1\ge0$$ The quadratic $\phi$ is convex because its leading coefficient is $\alpha-1\ge0$. Moreover, $\phi(1)=\alpha+\beta-4\ge0$ and $$\phi'(1)=2\alpha+\beta-4
=(\alpha+\beta-4)+\alpha\ge1$$ Hence $\phi(k)\ge0$ for every $k\ge1$, and therefore $w_{k+1}\le w_k \text{ } \forall \text{ } k\ge1.$ Similarly, the upper bound for $z_{k+1}/z_k$ is at most $1$ whenever, $$T(k)=(\alpha-1)k+(\beta-2)\ge0.$$ At $k=1$, $T(k)$ is $\alpha+\beta-3$ and as we assume $\alpha+\beta \geq4$ this means that $T(k)\geq1$, and $T(k)$ is non-decreasing in $k$ because $\alpha\ge1$. Consequently, $z_{k+1}\le z_k \text{ } \forall \text{ } k\ge1.$ Thus, both $\{w_k\}_{k\ge1}$ and $\{z_k\}_{k\ge1}$ are decreasing.

\smallskip
\noindent\textit{Exponential starlikeness.}
Because $\{c_{n}\}$ and $\{b_{n}\}$ are decreasing, for $z\in\mathcal{D}$,
\begin{align*}
\left|\frac{\mathbb{F}_{\alpha,\beta}^{\gamma}(z)}{z}\right|
  &\ge 1-\sum_{n=1}^{\infty}\frac{c_{n}}{n!}
   \ge 1-c_{1}(\alpha,\beta,\gamma)(e-1),\\
\left|{\mathbb{F}_{\alpha,\beta}^{\gamma}}'(z)
      -\frac{\mathbb{F}_{\alpha,\beta}^{\gamma}(z)}{z}\right|
  &\le \sum_{n=1}^{\infty}\frac{b_{n}}{n!}
   \le b_{1}(\alpha,\beta,\gamma)(e-1).
\end{align*}
Combining these and applying Lemma~\ref{Subordination Lemma}, we obtain
\begin{equation}\label{eq:star1}
\left|\frac{z{\mathbb{F}_{\alpha,\beta}^{\gamma}}'(z)}
        {\mathbb{F}_{\alpha,\beta}^{\gamma}(z)}-1\right|
  \le \frac{b_{1}(\alpha,\beta,\gamma)(e-1)}
         {1-c_{1}(\alpha,\beta,\gamma)(e-1)}
  < \frac{e-1}{e}.
\end{equation}
This gives $e\,b_{1}+c_{1}(e-1)<1$, and since
$b_{1}=c_{1}=[\Gamma(\beta)]^{\gamma}/[\Gamma(\alpha+\beta)]^{\gamma}$, we deduce
$$(2e-1)[\Gamma(\beta)]^{\gamma} < [\Gamma(\alpha+\beta)]^{\gamma}.$$

\smallskip
\noindent\textit{Exponential convexity.}
Since $\{d_{n}\}$ and $\{g_{n}\}$ are decreasing, for $z\in\mathcal{D}$,
\begin{align*}
|z{\mathbb{F}_{\alpha,\beta}^{\gamma}}''(z)|
  &\le \sum_{n=1}^{\infty}\frac{d_{n}}{(n-1)!}|z|^{n}
   < d_{1}(\alpha,\beta,\gamma)\,e,\\
|{\mathbb{F}_{\alpha,\beta}^{\gamma}}'(z)|
  &\ge 1-\sum_{n=1}^{\infty}\frac{g_{n}}{n!}
   \ge 1-g_{1}(\alpha,\beta,\gamma)(e-1).
\end{align*}
Therefore
\begin{equation}\label{eq:conv1}
\left|\frac{z{\mathbb{F}_{\alpha,\beta}^{\gamma}}''(z)}
        {{\mathbb{F}_{\alpha,\beta}^{\gamma}}'(z)}\right|
  < \frac{d_{1}(\alpha,\beta,\gamma)\,e}
         {1-g_{1}(\alpha,\beta,\gamma)(e-1)}
  < \frac{e-1}{e}.
\end{equation}
Since $d_{1}=g_{1}=2[\Gamma(\beta)]^{\gamma}/[\Gamma(\alpha+\beta)]^{\gamma}$,
this is equivalent to
$$(4e^{2}-4e+2)[\Gamma(\beta)]^{\gamma}
   < (e-1)[\Gamma(\alpha+\beta)]^{\gamma}. $$
\end{proof}
Thus, we have provided sufficient conditions for the exponential starlikeness and exponential convexity of the Le Roy type Mittag-Leffler functions and obtained sufficient coefficient bounds for the given geometric properties. In the next theorem, exponential starlikeness and convexity have been obtained for multi-index Le Roy type functions.

\begin{example}\label{ex:starlike-not-convex}
Let $\alpha=3$, $\beta=1$, and $\gamma=1$. This set of values satisfies all the prerequisite conditions of Theorem~\ref{thm:exp-starlike-convex} and the function is \textbf{exponentially starlike}, however the coefficient condition for exponential convexity fails because
$$(4e^{2}-4e+2)[\Gamma(1)]^{1}=4e^{2}-4e+2 > (e-1)[\Gamma(4)]^{1}=6(e-1).$$
In Figure~\ref{thm3}, the red region is the image of $\mathcal{D}$ under $z{\mathbb{F}_{\alpha,\beta}^{\gamma}}'(z) /\mathbb{F}_{\alpha,\beta}^{\gamma}(z) $ and the yellow region represents the image of $\mathcal{D}$ under $e^{z}$. The containment of the red region within the yellow illustrates exponential starlikeness. This demonstrates a case where the function meets the starlike criteria but not the convexity criteria.
\begin{figure}[h!]
    \centering
    \includegraphics[width=0.5\linewidth]{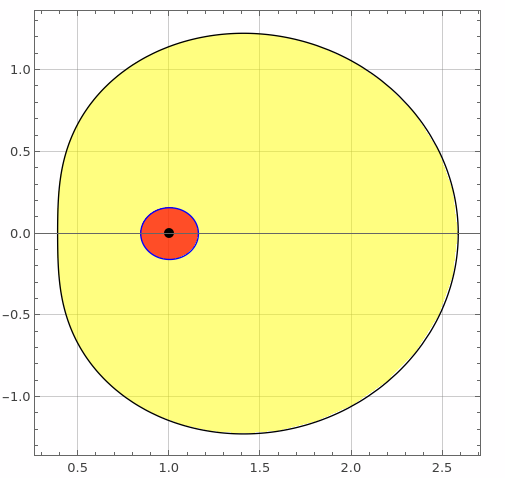}
    \caption{$\alpha=3$, $\beta=1$, $\gamma=1$: exponential starlikeness (red) inside $e^{z}$ (yellow)}
    \label{thm3}
\end{figure}
\end{example}
\begin{example}\label{ex:starlike-and-convex}
Let $\alpha=4$, $\beta=1$, and $\gamma=1$. This set of values satisfies all the prerequisite conditions of Theorem~\ref{thm:exp-starlike-convex}. Furthermore, these parameters satisfy the criteria for the function to be \textbf{exponentially starlike} and also satisfy the more restrictive condition for it to be \textbf{exponentially convex}. The exponential starlikeness is verified graphically in Figure~\ref{thm3all}, showing that the image of $z{\mathbb{F}_{4,1}^{1}}'(z)/\mathbb{F}_{4,1}^{1}(z)$ lies inside the image of $e^{z}$, the exponential convexity is verified in Figure~\ref{thm3all2}, showing that the image of $1+z{\mathbb{F}_{4,1}^{1}}''(z)/{\mathbb{F}_{4,1}^{1}}'(z)$ also lies inside the image of $e^{z}$. Thus the function belongs to both classes $S_{e}$ and $K_{e}$.
\begin{figure}[h!]
    \centering
    \begin{subfigure}{0.48\linewidth}
        \centering
        \includegraphics[width=\linewidth]{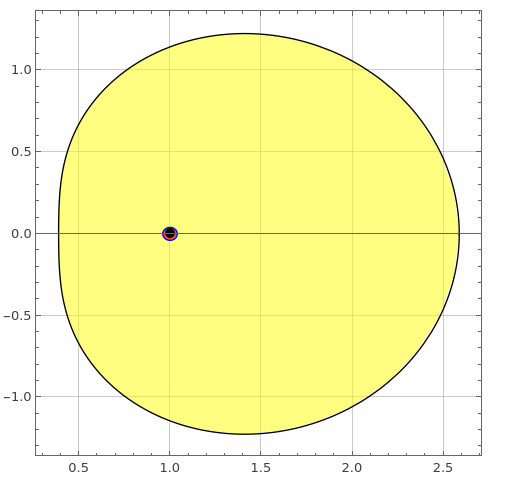}
        \caption{Exponential starlikeness: $z\frac{{\mathbb{F}_{4,1}^{1}}'(z)}{\mathbb{F}_{4,1}^{1}(z)}\subset e^{z}$}
        \label{thm3all}
    \end{subfigure}
    \hfill
    \begin{subfigure}{0.48\linewidth}
        \centering
        \includegraphics[width=\linewidth]{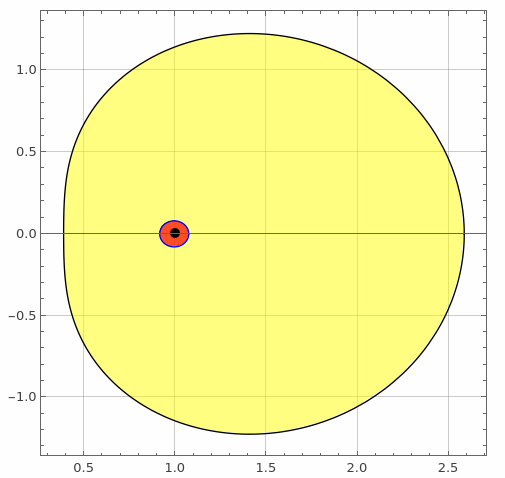}
        \caption{Exponential convexity: $1+z\frac{{\mathbb{F}_{4,1}^{1}}''(z)}{{\mathbb{F}_{4,1}^{1}}'(z)}\subset e^{z}$}
        \label{thm3all2}
    \end{subfigure}
    \caption{$\alpha=4$, $\beta=1$, $\gamma=1$: exponential starlikeness (left) and exponential convexity (right), in each case the red region is the relevant quantity subordinated to $e^{z}$ (yellow)}
\end{figure}
\end{example}

\begin{theorem}\label{thm:multivariate-exp}
Let $\mathbf{F}(z)=\mathbb{F}^{\gamma_{1},\gamma_{2},\ldots,\gamma_{n}}_{\alpha_{1},\alpha_{2},\ldots,\alpha_{n},\beta_{1},\beta_{2},\ldots,\beta_{n}}(z)$
be the standard normalized multi-index Le Roy type function. Let $\min(\alpha_i,\gamma_i)\ge 1$ and $\beta_i>0$ for all $i=1,\dots,n$,
and suppose $\alpha_j+\beta_j\ge 4$ for some $j\in\{1,\dots,n\}$. Then the 
following statements are true.
\begin{enumerate}
\item The normalized multi-index Le Roy type function $\mathbf{F}$ is exponentially starlike in 
$\mathcal{D}$ if it satisfies the coefficient condition
$$(2e-1){\prod_{i=1}^{n}\Gamma(\beta_{i})^{\gamma_{i}}} 
< {\prod_{i=1}^{n}[\Gamma(\alpha_{i}+\beta_{i})]^{\gamma_{i}}}.$$
\item The normalized multi-index Le Roy type function $\mathbf{F}$ is exponentially convex in 
$\mathcal{D}$ if it satisfies the coefficient condition
$$(4e^{2}-4e+2){\prod_{i=1}^{n}\Gamma(\beta_{i})^{\gamma_{i}}} 
< (e-1){\prod_{i=1}^{n}[\Gamma(\alpha_{i}+\beta_{i})]^{\gamma_{i}}}.$$
\end{enumerate}
\end{theorem}

\begin{proof}
Consider the normalized multivariate Le Roy type function
$$\mathbf{F}(z)=z+\sum_{k=2}^{\infty}
\frac{{\prod_{i=1}^{n}[\Gamma(\beta_{i})]^{\gamma_{i}}}z^{k}}
{\prod_{i=1}^{n}[\Gamma(\alpha_{i}(k-1)+\beta_{i})]^{\gamma_{i}}},$$
and let
$$x_{k} := \frac{\Gamma(k+1)}{\prod_{i=1}^{n}[\Gamma(\alpha_{i}k+\beta_{i})]^{\gamma_{i}}}, 
\qquad k\ge 1,$$
as in Lemma~\ref{newlem}. Set $w_k := kx_k$ and $z_k := (k+1)x_k$.\\[1em]
\noindent\textit{Monotonicity of $w_k$ and $z_k$.}
By the proof of Lemma~\ref{newlem}, under $\min(\alpha_i,\gamma_i)\ge 1$,
$$\frac{x_{k+1}}{x_{k}} \le \frac{k+1}{\alpha_jk+\beta_j}$$
for the index $j$ with $\alpha_j+\beta_j\ge 4$. Hence
$$\frac{w_{k+1}}{w_{k}} = \frac{k+1}{k}\cdot\frac{x_{k+1}}{x_{k}} 
\le \frac{(k+1)^{2}}{k(\alpha_jk+\beta_j)}.$$
This upper bound for $w_{k+1}/w_{k}$ is at most 1 whenever
$$\varphi(k):=(\alpha_j-1)k^{2}+(\beta_j-2)k-1 \ge 0.$$
Since $\alpha_j\ge 1$, $\varphi$ is convex, and
$$\varphi(1)=\alpha_j+\beta_j-4 \ge 0, \qquad 
\varphi'(1) = 2\alpha_j+\beta_j-4 = (\alpha_j+\beta_j -4)+\alpha_j \ge 1 > 0$$
by hypothesis. A convex function with $\varphi(1)\ge 0$ and $\varphi'(1)\ge 0$ 
satisfies $\varphi(k)\ge 0$ for all $k\ge 1$. Hence $w_{k+1}\le w_{k} \text{ } \forall \text{ } k \geq 1$, so 
$\{w_k\}_{k\geq1}$ is decreasing.\\

\noindent Similarly,
$$\frac{z_{k+1}}{z_{k}} = \frac{k+2}{k+1}\cdot\frac{x_{k+1}}{x_{k}} 
\le \frac{k+2}{\alpha_jk+\beta_j}$$ The upper bound on the  right hand side of the above inequality is at most $1$ whenever, $$T'(k)=(\alpha_j-1)k+(\beta_j-2)\ge0$$ So, at $k=1$, $T'(k) = \alpha_j+\beta_j-3\ge1$, and it is non-decreasing in $k$ as $\alpha_j\geq1$. Hence $z_{k+1}\leq z_{k} \text{ } \forall \text{ } k \geq 1$ which implies that
$\{z_k\}_{k\geq1}$ is also decreasing.

\smallskip
\noindent\textit{Part (1).} Since $c_n = x_n\prod_i\Gamma(\beta_i)^{\gamma_i}$ 
and $b_n = w_n\prod_i\Gamma(\beta_i)^{\gamma_i}$, both $\{c_n\}$ and 
$\{b_n\}$ are decreasing for $n\ge 1$ by Lemma~\ref{newlem} and the 
monotonicity of $\{w_k\}$ established above. Hence, for $z\in\mathcal{D}$,
$$\left|\frac{\mathbf{F}(z)}{z}\right| 
\ge 1-\sum_{n=1}^{\infty}\frac{c_{n}}{n!}
\ge 1-c_{1}(\alpha_{i},\beta_{i},\gamma_{i})(e-1),$$
$$\left|\mathbf{F}^{\prime}(z)-\frac{\mathbf{F}(z)}{z}\right| 
\le \sum_{n=1}^{\infty}\frac{b_{n}}{n!}
\le b_{1}(\alpha_{i},\beta_{i},\gamma_{i})(e-1).$$
Combining these and applying Lemma~\ref{Subordination Lemma}, we obtain
\begin{equation}\label{eq1}
\left|\frac{z\mathbf{F}^{\prime}(z)}{\mathbf{F}(z)}-1\right| 
\le \frac{b_{1}(\alpha_{i},\beta_{i},\gamma_{i})(e-1)}
{1-c_{1}(\alpha_{i},\beta_{i},\gamma_{i})(e-1)} < \frac{e-1}{e}.
\end{equation}
This gives $e\,b_{1}+c_{1}(e-1)<1$, and since 
$$b_{1}=c_{1}={\prod_{i=1}^{n}\Gamma(\beta_{i})^{\gamma_{i}}}
/{\prod_{i=1}^{n}[\Gamma(\alpha_{i}+\beta_{i})]^{\gamma_{i}}}$$ we deduce that
$$(2e-1){\prod_{i=1}^{n}\Gamma(\beta_{i})^{\gamma_{i}}} 
< {\prod_{i=1}^{n}[\Gamma(\alpha_{i}+\beta_{i})]^{\gamma_{i}}}.$$

\smallskip
\noindent\textit{Part (2).} Since $d_n = g_n = z_n\prod_i\Gamma(\beta_i)^{\gamma_i}$,
both $\{d_n\}$ and $\{g_n\}$ are decreasing for $n\ge 1$ by the 
monotonicity of $\{z_k\}$ established above. Hence, for $z\in\mathcal{D}$,
$$|z\mathbf{F}^{\prime\prime}(z)| 
\le \sum_{n=1}^{\infty}\frac{d_{n}}{(n-1)!}|z|^{n} 
< d_{1}(\alpha_{i},\beta_{i},\gamma_{i})\,e,$$
$$|\mathbf{F}^{\prime}(z)| 
\ge 1-\sum_{n=1}^{\infty}\frac{g_{n}}{n!} 
\ge 1-g_{1}(\alpha_{i},\beta_{i},\gamma_{i})(e-1).$$
Therefore
\begin{equation}\label{eq2}
\left|\frac{z\mathbf{F}^{\prime\prime}(z)}{\mathbf{F}^{\prime}(z)}\right| 
< \frac{d_{1}(\alpha_{i},\beta_{i},\gamma_{i})\,e}
{1-g_{1}(\alpha_{i},\beta_{i},\gamma_{i})(e-1)} < \frac{e-1}{e}.
\end{equation}
Since $d_{1}=g_{1}=2{\prod_{i=1}^{n}\Gamma(\beta_{i})^{\gamma_{i}}}
/{\prod_{i=1}^{n}[\Gamma(\alpha_{i}+\beta_{i})]^{\gamma_{i}}}$, this is 
equivalent to
$$(4e^{2}-4e+2){\prod_{i=1}^{n}\Gamma(\beta_{i})^{\gamma_{i}}} 
< (e-1){\prod_{i=1}^{n}[\Gamma(\alpha_{i}+\beta_{i})]^{\gamma_{i}}}. 
$$
\end{proof}

\begin{example}
Let $n=3$ and consider the following parameter sets:
$(\alpha_1, \beta_1, \gamma_1) = (3, 1, 1)$,
$(\alpha_2, \beta_2, \gamma_2) = (2, 1, 1)$,
$(\alpha_3, \beta_3, \gamma_3) = (2, 1, 1)$.
Here $\min(\alpha_i,\gamma_i)\ge1$, $\beta_i>0$,
$\alpha_1+\beta_1=4\ge4$, and
$\prod_{i=1}^{3}\Gamma(\beta_i)^{\gamma_i}=1$.
Since
$\prod_{i=1}^{3}[\Gamma(\alpha_i+\beta_i)]^{\gamma_i}=
\Gamma(4)\cdot\Gamma(3)\cdot\Gamma(3)=24$,
both the starlikeness condition $(2e-1)<24$ and the convexity condition
$(4e^2-4e+2)<(e-1)\cdot24$ hold.
Thus all the hypotheses of Theorem~\ref{thm:multivariate-exp} are satisfied.
The exponential starlikeness is verified graphically in Figure~\ref{thm4},
showing that the image of $z\mathbf{F}'(z)/\mathbf{F}(z)$ lies inside $e^{z}$, the exponential convexity is verified in Figure~\ref{thm4all}, showing that
$1+z\mathbf{F}''(z)/\mathbf{F}'(z)$ also lies inside $e^{z}$.

\begin{figure}[h!]
    \centering
    \begin{subfigure}{0.48\linewidth}
        \centering
        \includegraphics[width=\linewidth]{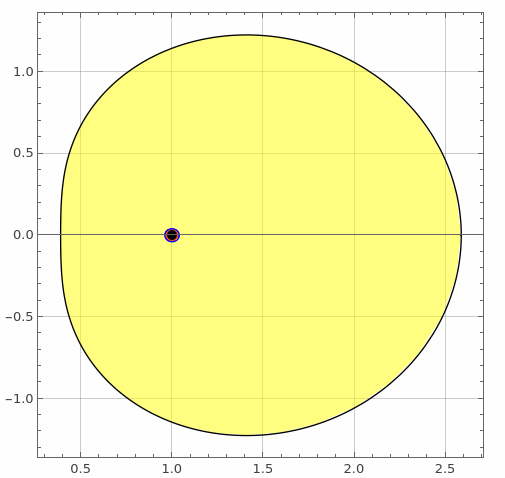}
        \caption{Exponential starlikeness: $z\frac{\mathbf{F}'(z)}{\mathbf{F}(z)}\subset e^{z}$}
        \label{thm4}
    \end{subfigure}
    \hfill
    \begin{subfigure}{0.48\linewidth}
        \centering
        \includegraphics[width=\linewidth]{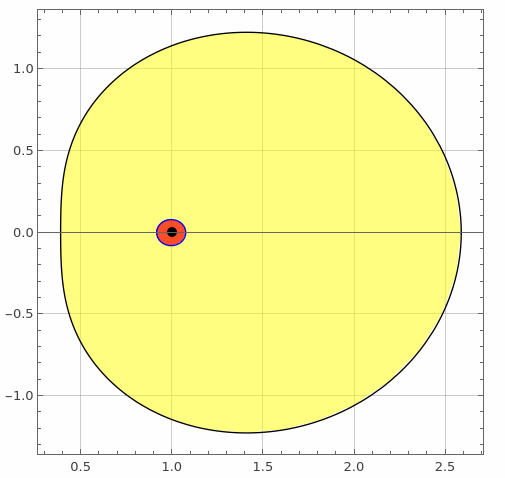}
        \caption{Exponential convexity:  $1+z\frac{\mathbf{F}''(z)}{\mathbf{F}'(z)}\subset e^{z}$}
        \label{thm4all}
    \end{subfigure}
    \caption{$n=3$ with $(\alpha_1,\beta_1,\gamma_1)=(3,1,1)$,
    $(\alpha_2,\beta_2,\gamma_2)=(2,1,1)$,
    $(\alpha_3,\beta_3,\gamma_3)=(2,1,1)$:
    exponential starlikeness (left) and convexity (right), red region inside $e^{z}$ (yellow)}
\end{figure}
\end{example}

\begin{corollary}\label{cor:multivariate-starlike}
Let $\min(\alpha_i,\gamma_i)\ge 1$ and $\beta_i>0$ for all $i=1,\dots,n$,
and suppose $\alpha_j+\beta_j\ge 4$ for some $j\in\{1,\dots,n\}$.
Then $\mathbf{F}$ is starlike in $\mathcal{D}$ if
$$
2(e-1)\prod_{i=1}^{n}\Gamma(\beta_{i})^{\gamma_{i}}
< \prod_{i=1}^{n}[\Gamma(\alpha_{i}+\beta_{i})]^{\gamma_{i}}.
$$
\end{corollary}

\begin{proof}
Using definition of starlikeness, it suffices to establish 
$\left|\frac{z\mathbf{F}'(z)}{\mathbf{F}(z)}-1\right|<1$. From \eqref{eq1}, this holds whenever
$$
\frac{b_{1}(e-1)}{1-c_{1}(e-1)} < 1.
$$
Since $1 - c_1(e-1) > 0$ under the given condition, cross-multiplying yields
$$
(b_{1}+c_{1})(e-1) < 1.
$$
Substituting
$b_{1}=c_{1}={\prod_{i=1}^{n}\Gamma(\beta_{i})^{\gamma_{i}}}\big/
{\prod_{i=1}^{n}[\Gamma(\alpha_{i}+\beta_{i})]^{\gamma_{i}}}$
gives
$$
2(e-1)\prod_{i=1}^{n}\Gamma(\beta_{i})^{\gamma_{i}}
< \prod_{i=1}^{n}[\Gamma(\alpha_{i}+\beta_{i})]^{\gamma_{i}}.
$$
\end{proof}

\begin{example}
Let $n=3$ and consider the parameter sets
$(\alpha_1,\beta_1,\gamma_1)=(3.0,\,1.0,\,1.2)$,
$(\alpha_2,\beta_2,\gamma_2)=(1.5,\,1.5,\,1.0)$,
$(\alpha_3,\beta_3,\gamma_3)=(2.0,\,1.0,\,1.0)$.
Here $\min(\alpha_i,\gamma_i)\ge1$, $\beta_i>0$, $\alpha_1+\beta_1=4\ge4$, and
$\prod_{i=1}^{3}\Gamma(\beta_i)^{\gamma_i}=\Gamma(1.5)\approx0.886$.
Since
$\prod_{i=1}^{3}[\Gamma(\alpha_i+\beta_i)]^{\gamma_i}
 =\Gamma(4)^{1.2}\cdot\Gamma(3)\cdot\Gamma(3)\approx34.34$,
the condition $2(e-1)\prod\Gamma(\beta_i)^{\gamma_i}<\prod[\Gamma(\alpha_i+\beta_i)]^{\gamma_i}$
is satisfied, so by Corollary~\ref{cor:multivariate-starlike} the function is starlike in $\mathcal{D}$. In Figure~\ref{cor_starlike}, the orange region is the image of the unit disk $\mathcal{D}$ under the normalized function $\mathbf{F}(z)$, with the boundary shown in blue. The dotted unit circle is shown for reference to compare with the next example. The image is a simple closed curve encircling the origin, confirming that $\mathbf{F}$ maps $\mathcal{D}$ onto a starlike domain.
\begin{figure}[h!]
    \centering
    \includegraphics[width=0.5\linewidth]{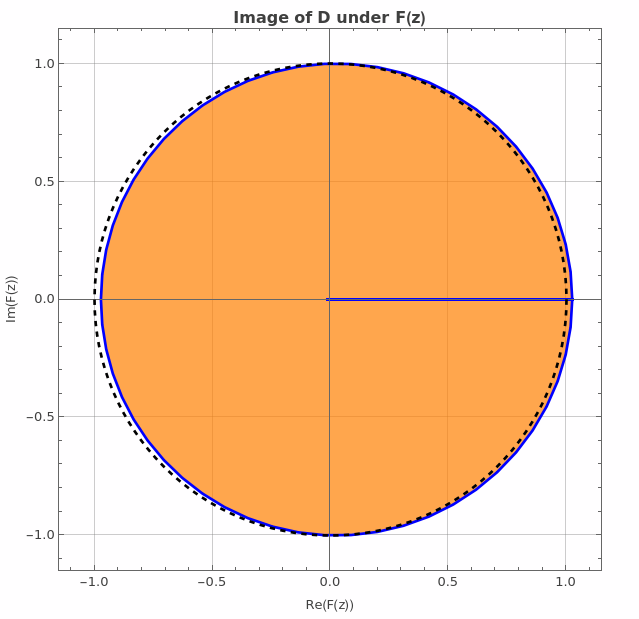}
    \caption{$n=3$ with $(\alpha_1,\beta_1,\gamma_1)=(3.0,1.0,1.2)$,
    $(\alpha_2,\beta_2,\gamma_2)=(1.5,1.5,1.0)$,
    $(\alpha_3,\beta_3,\gamma_3)=(2.0,1.0,1.0)$: image of $\mathcal{D}$ under $\mathbf{F}(z)$ (orange) with boundary (blue), the dotted circle is the unit circle}
    \label{cor_starlike}
\end{figure}
\end{example}

\begin{corollary}\label{cor:multivariate-convex}
Let $\min(\alpha_i,\gamma_i)\ge 1$ and $\beta_i>0$ for all $i=1,\dots,n$,
and suppose $\alpha_j+\beta_j\ge 4$ for some $j\in\{1,\dots,n\}$.
Then $\mathbf{F}$ is convex in $\mathcal{D}$ if
$$
2(2e-1)\prod_{i=1}^{n}\Gamma(\beta_{i})^{\gamma_{i}}
< \prod_{i=1}^{n}[\Gamma(\alpha_{i}+\beta_{i})]^{\gamma_{i}}.
$$
\end{corollary}

\begin{proof}
In the similar manner by definition of convexity, it suffices to obtain
$\left|\frac{z\mathbf{F}''(z)}{\mathbf{F}'(z)}\right| < 1$.
From \eqref{eq2}, this holds whenever
$$
\frac{d_{1}\,e}{1-g_{1}(e-1)} < 1.
$$
Since $1-g_{1}(e-1)>0$ under the stated condition, this is equivalent to
$$
d_{1}\,e + g_{1}(e-1) < 1.
$$
Substituting $d_{1}=g_{1}=
2\prod_{i=1}^{n}\Gamma(\beta_{i})^{\gamma_{i}}\big/
\prod_{i=1}^{n}[\Gamma(\alpha_{i}+\beta_{i})]^{\gamma_{i}}$
and factoring gives
$$
d_{1}\left(e+(e-1)\right) < 1,
$$
which gives
$$
2(2e-1)\prod_{i=1}^{n}\Gamma(\beta_{i})^{\gamma_{i}}
< \prod_{i=1}^{n}[\Gamma(\alpha_{i}+\beta_{i})]^{\gamma_{i}}.
$$
\end{proof}

\begin{example}
Let $n=3$ and consider the parameter sets
$(\alpha_1,\beta_1,\gamma_1)=(2.5,\,1.0,\,1.0)$,
$(\alpha_2,\beta_2,\gamma_2)=(2.0,\,1.0,\,1.5)$,
$(\alpha_3,\beta_3,\gamma_3)=(3.0,\,1.0,\,1.0)$.
Here $\min(\alpha_i,\gamma_i)\ge1$, $\beta_i>0$, $\alpha_3+\beta_3=4\ge4$, and
$\prod_{i=1}^{3}\Gamma(\beta_i)^{\gamma_i}=1$.
Since
$\prod_{i=1}^{3}[\Gamma(\alpha_i+\beta_i)]^{\gamma_i}
 =\Gamma(3.5)\cdot\Gamma(3)^{1.5}\cdot\Gamma(4)\approx56.40$,
the condition $2(2e-1)\prod\Gamma(\beta_i)^{\gamma_i}<\prod[\Gamma(\alpha_i+\beta_i)]^{\gamma_i}$
is satisfied, so by Corollary~\ref{cor:multivariate-convex} the function is convex in $\mathcal{D}$. In Figure~\ref{cor_convex}, the orange region is the image of the unit disk $\mathcal{D}$ under the normalized function $\mathbf{F}(z)$, with the boundary shown in blue and the dotted unit circle for reference. The dotted unit circle was graphed with the image of $\mathcal{D}$ under $\mathbf{F}(z)$ here and also in Figure~\ref{cor_starlike}, so the difference is visible even though it is very small.
\begin{figure}[h!]
    \centering
    \includegraphics[width=0.5\linewidth]{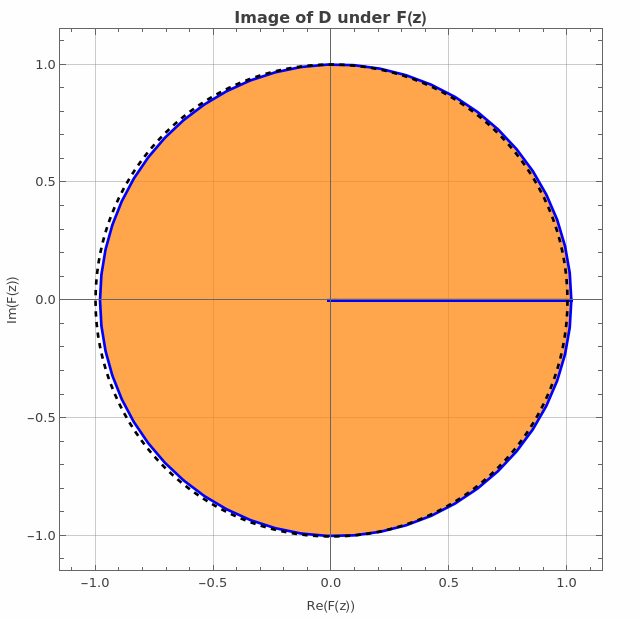}
    \caption{$n=3$ with $(\alpha_1,\beta_1,\gamma_1)=(2.5,1.0,1.0)$,
    $(\alpha_2,\beta_2,\gamma_2)=(2.0,1.0,1.5)$,
    $(\alpha_3,\beta_3,\gamma_3)=(3.0,1.0,1.0)$: image of $\mathcal{D}$ under $\mathbf{F}(z)$ (orange), the dotted circle is the unit circle}
    \label{cor_convex}
\end{figure}
\end{example}

\begin{corollary}\label{cor:multivariate-starlike-sharp}
Let $\min(\alpha_i,\gamma_i)\ge 1$ and $\beta_i>0$ for all $i=1,\dots,n$, and suppose $\alpha_j+\beta_j\ge 3$ for some $j\in\{1,\dots,n\}$. Then the normalized multi-index Le Roy type function $\mathbf{F}$ is starlike in 
$\mathcal{D}$ if
$$\sqrt{5}\,(e-1)\prod_{i=1}^{n}\Gamma(\beta_{i})^{\gamma_{i}}
< \prod_{i=1}^{n}[\Gamma(\alpha_{i}+\beta_{i})]^{\gamma_{i}}.$$
\end{corollary}

\begin{proof}
Consider the normalized multivariate Le Roy type function
$$\mathbf{F}(z) = z + \sum_{k=2}^{\infty}
\frac{{\prod_{i=1}^{n}[\Gamma(\beta_{i})]^{\gamma_{i}}}\,z^{k}}
{\prod_{i=1}^{n}[\Gamma(\alpha_{i}(k-1)+\beta_{i})]^{\gamma_{i}}}.$$
To establish starlikeness in $\mathcal{D}$, we apply Lemma~\ref{Starlike},
which requires $|\mathbf{F}'(z)-1| < \dfrac{2}{\sqrt{5}}$.
Substituting $m=k-1$ we obtain
$$\left|\mathbf{F}'(z)-1\right|
\le \sum_{m=1}^{\infty}\frac{g_{m}}{m!}|z|^{m},$$
where
$$g_{m}
=\frac{(m+1)\Gamma(m+1){\prod_{i=1}^{n}[\Gamma(\beta_{i})]^{\gamma_{i}}}}
{\prod_{i=1}^{n}[\Gamma(\alpha_{i}m+\beta_{i})]^{\gamma_{i}}}.$$
Let $x_{m}:=\Gamma(m+1)/\prod_{i=1}^{n}[\Gamma(\alpha_{i}m+\beta_{i})]^{\gamma_{i}}$
and $z_{m}:=(m+1)x_{m}$ as in the proof of Theorem~\ref{thm:multivariate-exp},
so that $g_{m}=z_{m}\prod_{i=1}^{n}[\Gamma(\beta_{i})]^{\gamma_{i}}$.
By the same ratio estimate as in Theorem~\ref{thm:multivariate-exp},
$\{z_{m}\}_{m\ge 1}$ is decreasing under $\min(\alpha_i,\gamma_i)\ge1$
and $\alpha_j+\beta_j\ge3$. Hence, replacing each coefficient by its
supremum $g_1$, we obtain for $z\in\mathcal{D}$,
$$\left|\mathbf{F}'(z)-1\right|
\le g_{1}\sum_{m=1}^{\infty}\frac{1}{m!}
= g_{1}(e-1).$$
Since
$g_{1}
=\dfrac{2{\prod_{i=1}^{n}\Gamma(\beta_{i})^{\gamma_{i}}}}
{\prod_{i=1}^{n}[\Gamma(\alpha_{i}+\beta_{i})]^{\gamma_{i}}},$
the condition $g_1(e-1) < \frac{2}{\sqrt{5}}$ is equivalent to
$$\sqrt{5}\,(e-1)\prod_{i=1}^{n}\Gamma(\beta_{i})^{\gamma_{i}}
< \prod_{i=1}^{n}[\Gamma(\alpha_{i}+\beta_{i})]^{\gamma_{i}}. $$
\end{proof}

\noindent Since $\sqrt{5}(e-1) > 2(e-1)$, the coefficient condition of 
Corollary~\ref{cor:multivariate-starlike-sharp} is strictly stronger than the 
starlikeness condition of Corollary~\ref{cor:multivariate-starlike}. The difference is due to the underlying sufficient conditions used. 
Corollary~\ref{cor:multivariate-starlike-sharp} invokes the sharper criterion 
$|\mathbf{F}'(z)-1|<\frac{2}{\sqrt{5}}$ of Lemma~\ref{Starlike}, whereas 
Corollary~\ref{cor:multivariate-starlike} relies on the elementary 
disk-containment criterion i.e.\ 
$\left|\frac{z\mathbf{F}'(z)}{\mathbf{F}(z)}-1\right|<1$.

On the other hand, the parameter hypothesis $\alpha_j+\beta_j\ge3$ of
Corollary~\ref{cor:multivariate-starlike-sharp} is weaker than the
$\alpha_j+\beta_j\ge4$ required in Theorem~\ref{thm:multivariate-exp}.
This is because the corollary uses only the decreasing property of
$\{z_{m}\}$ (which needs $\alpha_j+\beta_j\ge3$), whereas Theorem~\ref{thm:multivariate-exp}
must also guarantee that $\{w_{m}\}$ decreases, imposing the stricter
$\alpha_j+\beta_j\ge4$ condition.

\begin{example}
Let $n=3$ and consider the parameter sets
$(\alpha_1,\beta_1,\gamma_1)=(2.0,\,1.2,\,1.0)$,
$(\alpha_2,\beta_2,\gamma_2)=(1.5,\,1.0,\,1.5)$,
$(\alpha_3,\beta_3,\gamma_3)=(2.2,\,1.0,\,1.0)$.
Here $\min(\alpha_i,\gamma_i)\ge1$, $\beta_i>0$, $\alpha_3+\beta_3=3.2\ge3$, and
$\prod_{i=1}^{3}\Gamma(\beta_i)^{\gamma_i}=\Gamma(1.2)\approx0.918$.
Since
$\prod_{i=1}^{3}[\Gamma(\alpha_i+\beta_i)]^{\gamma_i}
 =\Gamma(3.2)\cdot\Gamma(2.5)^{1.5}\cdot\Gamma(3.2)\approx9.005$,
the condition $\sqrt{5}(e-1)\prod\Gamma(\beta_i)^{\gamma_i}<\prod[\Gamma(\alpha_i+\beta_i)]^{\gamma_i}$
is satisfied, so by Corollary~\ref{cor:multivariate-starlike-sharp} the function is starlike in $\mathcal{D}$. In Figure~\ref{star_thm}, the orange region is the image of the unit disk $\mathcal{D}$ under the normalized function $\mathbf{F}(z)$, with the boundary shown in blue. The dotted unit circle is again used as a reference to differentiate between the previous examples. The image  illustrates that $\mathbf{F}$ maps $\mathcal{D}$ onto a starlike domain.
\begin{figure}[h!]
    \centering
    \includegraphics[width=0.5\linewidth]{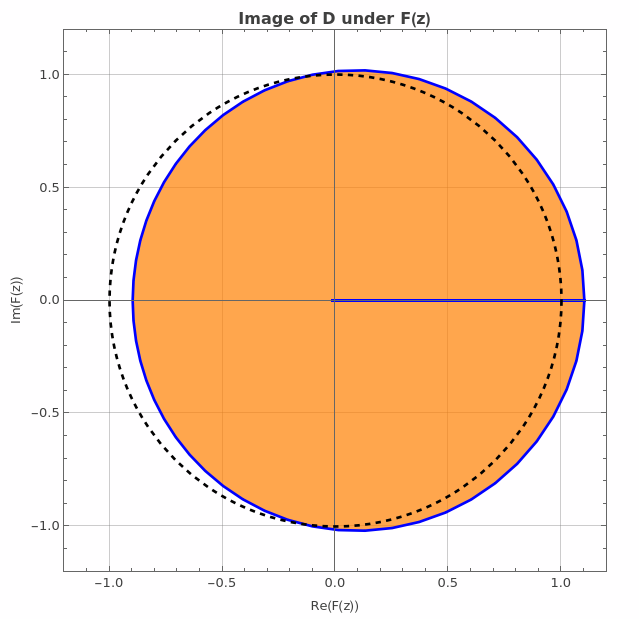}
    \caption{$n=3$ with $(\alpha_1,\beta_1,\gamma_1)=(2.0,1.2,1.0)$,
    $(\alpha_2,\beta_2,\gamma_2)=(1.5,1.0,1.5)$,
    $(\alpha_3,\beta_3,\gamma_3)=(2.2,1.0,1.0)$: image of $\mathcal{D}$ under $\mathbf{F}(z)$ (orange) with boundary (blue), the dotted circle is the unit circle}
    \label{star_thm}
\end{figure}
\end{example}

 \section{Close-to-convexity, Starlikeness and convexity  of Multivariate Le Roy Type Functions}
 In this section, geometric properties such as Close-to-convexity, starlikeness and convexity have been investigated for Multivariate Le Roy type functions.
We continue to denote the original multi-index function by
$$\mathcal{F}(z)=F^{\gamma_{1},\gamma_{2},\ldots,\gamma_{n}}_{\alpha_{1},\alpha_{2},\ldots,\alpha_{n},\beta_{1},\beta_{2},\ldots,\beta_{n}}(z)
=\sum_{m=0}^{\infty}\frac{z^m}{\prod_{i=1}^{n}[\Gamma(\alpha_i m+\beta_i)]^{\gamma_i}}$$
Its standard normalized form is
$$\mathbf{F}(z)=\mathbb{F}^{\gamma_{1},\gamma_{2},\ldots,\gamma_{n}}_{\alpha_{1},\alpha_{2},\ldots,\alpha_{n},\beta_{1},\beta_{2},\ldots,\beta_{n}}(z)
=z\prod_{i=1}^{n}[\Gamma(\beta_i)]^{\gamma_i}
\mathcal{F}(z)
=z+\sum_{m=2}^{\infty}
\frac{\prod_{i=1}^{n}[\Gamma(\beta_i)]^{\gamma_i}z^m}
{\prod_{i=1}^{n}[\Gamma(\alpha_i(m-1)+\beta_i)]^{\gamma_i}}.
$$
Consequently,
\begin{align*}
\mathbf{F}'(z)
&=1+\sum_{m=2}^{\infty}
\frac{m\prod_{i=1}^{n}[\Gamma(\beta_i)]^{\gamma_i}z^{m-1}}
{\prod_{i=1}^{n}[\Gamma(\alpha_i(m-1)+\beta_i)]^{\gamma_i}},\\
\mathbf{F}''(z)
&=\sum_{m=2}^{\infty}
\frac{m(m-1)\prod_{i=1}^{n}[\Gamma(\beta_i)]^{\gamma_i}z^{m-2}}
{\prod_{i=1}^{n}[\Gamma(\alpha_i(m-1)+\beta_i)]^{\gamma_i}}.
\end{align*}
The close-to-convexity result uses a different normalization which is derived in the theorem itself. Therefore, we denote it separately by
$$\mathbb{N}(z)
:=\prod_{i=1}^{n}[\Gamma(\alpha_i+\beta_i)]^{\gamma_i}
\left[\mathcal{F}(z)
-\frac{1}{\prod_{i=1}^{n}[\Gamma(\beta_i)]^{\gamma_i}}\right].$$

\begin{theorem}\label{thm:close-to-convex}
    Let $\min(\alpha_i,\gamma_i)\ge 1$ and $\beta_i>0$ for all $i=1,\dots,n$, and suppose $\alpha_j+\beta_j\ge 2$ for some $j\in\{1,\dots,n\}$. Then the normalized function $\mathbb{N}(z)$ is close-to-convex with respect to the starlike function $-\log (1-z)$ in $\mathcal{D}$,
    and consequently, it is univalent in $\mathcal{D}$.
\end{theorem}
\begin{proof}
Consider
$$
\mathcal{F}(z)=\sum_{m=0}^{\infty} \frac{z^{m}}{\prod_{i=1}^{n}[\Gamma(\alpha_{i}{m}+\beta_{i})]^{\gamma_{i}}}.
$$
Then
\begin{align*}
& \mathcal{F}(z)-\frac{1}{\prod_{i=1}^{n}[\Gamma(\beta_{i})]^{\gamma_{i}}}
   =\frac{z}{\prod_{i=1}^{n}[\Gamma(\alpha_{i}+\beta_{i})]^{\gamma_{i}}}
   +\frac{z^{2}}{\prod_{i=1}^{n}[\Gamma(2\alpha_{i}+\beta_{i})]^{\gamma_{i}}}+\cdots \\
   \text{implies}\\
& \prod_{i=1}^{n}[\Gamma(\alpha_{i}+\beta_{i})]^{\gamma_{i}}
   \left[\mathcal{F}(z)-\dfrac{1}{\prod_{i=1}^{n}[\Gamma(\beta_i)]^{\gamma_i}}\right] = z+\sum_{m=2}^{\infty}
   \frac{\prod_{i=1}^{n}[\Gamma(\alpha_{i}+\beta_{i})]^{\gamma_{i}}\,z^{m}}
        {\prod_{i=1}^{n}[\Gamma(\alpha_{i}m+\beta_{i})]^{\gamma_{i}}}
   =:\mathbb{N}(z).
\end{align*}
Set $a_{k}:=\prod_{i=1}^{n}[\Gamma(\alpha_{i}+\beta_{i})]^{\gamma_{i}}/
\prod_{i=1}^{n}[\Gamma(\alpha_{i}k+\beta_{i})]^{\gamma_{i}}$
for $k\ge 2$, and $a_{1}=1$.
By Ozaki's Lemma it suffices to verify
\begin{equation}\label{eq:ozaki}
1\ge 2a_{2}\ge 3a_{3}\ge\cdots\ge k a_{k}\ge (k+1)a_{k+1}\ge\cdots\ge 0.
\end{equation}
For the chain $k a_{k}\ge (k+1)a_{k+1}$, we use a ratio estimate.
For each $i$ and every $k\ge1$, let $x:=\alpha_i k+\beta_i$. Then $x+1\ge 2$ and, since $\Gamma$ is increasing on $[2,\infty)$ and $\alpha_i\ge1$, we have $x+\alpha_i\ge x+1$, hence $\Gamma(\alpha_i(k+1)+\beta_i)=\Gamma(x+\alpha_i)
   \ge \Gamma(x+1)
   =x\,\Gamma(x).$
Consequently,
\begin{align*}
\frac{k a_{k}}{(k+1)a_{k+1}}
   &= \frac{k}{k+1}\,
      \prod_{i=1}^{n}\frac{[\Gamma(\alpha_i(k+1)+\beta_i)]^{\gamma_i}}
                       {[\Gamma(\alpha_i k+\beta_i)]^{\gamma_i}}\ge \frac{k}{k+1}\,\prod_{i=1}^{n}x^{\gamma_i}
    = \frac{k}{k+1}\,\prod_{i=1}^{n}(\alpha_i k+\beta_i)^{\gamma_i}.
\end{align*}
Since $\alpha_i k+\beta_i\ge 1$ and $\gamma_i\ge1$,
$(\alpha_i k+\beta_i)^{\gamma_i}\ge\alpha_i k+\beta_i$, and therefore
$$\prod_{i=1}^{n}(\alpha_i k+\beta_i)^{\gamma_i}
   \ge \alpha_j k+\beta_j$$
for the index $j$ with $\alpha_j+\beta_j\ge2$.
Thus
$$\frac{k a_{k}}{(k+1)a_{k+1}}
   \ge \frac{k}{k+1}\,(\alpha_j k+\beta_j).$$
The right-hand side is $\ge1$ whenever
$\alpha_j k+\beta_j\ge\frac{k+1}{k}=1+\frac1k$,
which holds for all $k\ge1$ because, at $k=1$ the hypothesis
$\alpha_j+\beta_j\ge2$ gives the required bound, and for $k\ge2$
we have $\alpha_j k+\beta_j\ge k>1+\frac1k$ (since $\alpha_j\ge1$ and $\beta_j>0$).
Hence $k a_{k}\ge(k+1)a_{k+1}$ for every $k\ge1$, and in particular
$1\ge 2a_{2}$ (because $a_{1}=1$).
Thus the entire chain \eqref{eq:ozaki} is satisfied.

By Stirling's formula,
$\Gamma(\alpha_i k+\beta_i)\sim\sqrt{2\pi}(\alpha_i k)^{\alpha_i k+\beta_i-1/2}e^{-\alpha_i k}$
as $k\to\infty$, so each factor $[\Gamma(\alpha_i k+\beta_i)]^{\gamma_i}$ grows
super-exponentially. Thus, $k a_{k}\to0$ as $k\to\infty$.
So the condition \ref{eq:ozaki} is satisfied,
hence $\mathbb{N}(z)$ is close-to-convex with respect to
$-\log(1-z)$ in $\mathcal{D}$.
\end{proof}

\begin{example}
Let $n=3$ and consider the parameter sets
$(\alpha_1,\beta_1,\gamma_1)=(1.1,\,1.3,\,1.6)$,
$(\alpha_2,\beta_2,\gamma_2)=(1.2,\,1.1,\,1.05)$,
$(\alpha_3,\beta_3,\gamma_3)=(1.1,\,1.2,\,1.05)$.
Here $\min(\alpha_i,\gamma_i)\ge1$, $\beta_i>0$, and $\alpha_2+\beta_2=2.3\ge2$.
Thus all hypotheses of Theorem~\ref{thm:close-to-convex} are satisfied, so the function is close-to-convex and therefore univalent in $\mathcal{D}$. Figure~\ref{close_convex} shows the image of the unit disk $\mathcal{D}$ under $\mathbb{N}(z)$, the curve winds around the origin without self-intersection, illustrating the close-to-convex property.
\begin{figure}[h!]
    \centering
    \includegraphics[width=0.5\linewidth]{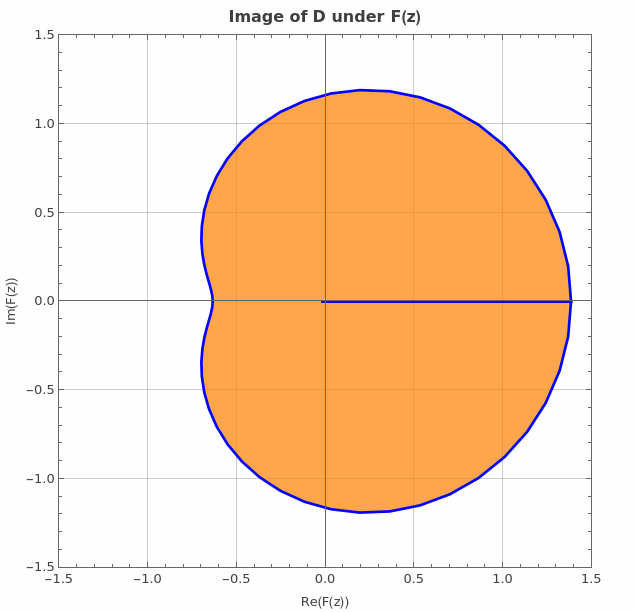}
    \caption{$n=3$ with $(\alpha_1,\beta_1,\gamma_1)=(1.1,1.3,1.6)$,
    $(\alpha_2,\beta_2,\gamma_2)=(1.2,1.1,1.05)$,
    $(\alpha_3,\beta_3,\gamma_3)=(1.1,1.2,1.05)$:
    image of the unit disk under $\mathbb{N}(z)$ (orange), illustrating close-to-convexity}
    \label{close_convex}
\end{figure}
\end{example}

\begin{theorem}
    Let $\min (\alpha_i, \gamma_i) \geq 1, \beta_i>0$ for all $i=1,2,\cdots,n$ such that $\alpha_j+\beta_j \geq 2$ for some $j \in {1,2,\cdots,n}$. Then the following assertions hold true,
    \begin{enumerate}[label=(\alph*)]
        \item If $(e-1)\prod_{i=1}^n[\Gamma(\beta_i)]^{\gamma_i}<\prod_{i=1}^n[\Gamma(\alpha_i+\beta_i)]^{\gamma_i}$, then the function $\mathbf{F}(z)$ is starlike in $\mathcal{D}_{1 / 2}$.
        \item If $2(e-1)\prod_{i=1}^n[\Gamma(\beta_i)]^{\gamma_i}<\prod_{i=1}^n[\Gamma(\alpha_i+\beta_i)]^{\gamma_i}$ and $\beta_i \geq 2$ for all $i =1,2,\cdots,n$, then the function $\mathbf{F}(z)$ is convex in $\mathcal{D}_{1 / 2}$.
    \end{enumerate}
\end{theorem}
\begin{proof}
    \begin{enumerate}[label=(\alph*)]
        \item Using Lemma \ref{newlem} we calculate,

\begin{align*}
\left|\frac{\mathbf{F}(z)}{z}-1\right| & \leq \sum_{k=1}^{\infty} A_{k+1}|z|^{k} \\
& \leq \mathbf{F}(1)-1  \\
& \leq \theta^{\gamma_{1},\gamma_{2},\ldots,\gamma_{n}}_{\alpha_{1},\alpha_{2},\ldots,\alpha_{n},\beta_{1},\beta_{2},\ldots, \beta_{n}}(e-1),
\end{align*}
for all $z \in \mathcal{D}$, where $A_{k+1}=\prod_{i=1}^{n}\left[\frac{\Gamma(\beta_i)}{\Gamma(\alpha_i k+\beta_i)}\right]^{\gamma_i}$. Hence, under the given hypotheses we obtain $
\left|\frac{\mathbf{F}(z)}{z}-1\right|<1, z \in \mathcal{D}
$ and consequently the function $\mathbf{F}(z)$ is starlike in $\mathcal{D}_{1 / 2}$ by the means of Lemma \ref{mcgregorlemma}.\\
\item We compute the following,
\begin{align*}
\left(\mathbf{F}(z)\right)^{\prime}-1 & =\sum_{k=1}^{\infty}(k+1) A_{k+1} z^{k} =\sum_{k=1}^{\infty} \frac{y_{k} z^{k}}{k!}
\end{align*}
where 
\begin{equation*}
y_{k}=\frac{\prod_{i=1}^n[\Gamma(\beta_i)]^{\gamma_i} \Gamma(k+2)}{\prod_{i=1}^n[\Gamma(\alpha_i k+\beta_i)]^{\gamma_i}}, k \geq 1 
\end{equation*}
We now define the function, $$f^{\gamma_{1},\gamma_{2},\ldots,\gamma_{n}}_{\alpha_{1},\alpha_{2},\ldots,\alpha_{n},\beta_{1},\beta_{2},\ldots, \beta_{n}}(x)=\frac{\Gamma(x+2)}{\prod_{i=1}^n[\Gamma(\alpha_i x+\beta_i)]^{\gamma_i}}, x>0 $$ 
Therefore,
\begin{equation*}
\left(f^{\gamma_{1},\gamma_{2},\ldots,\gamma_{n}}_{\alpha_{1},\alpha_{2},\ldots,\alpha_{n},\beta_{1},\beta_{2},\ldots, \beta_{n}}(x)\right)^{\prime}=f^{\gamma_{1},\gamma_{2},\ldots,\gamma_{n}}_{\alpha_{1},\alpha_{2},\ldots,\alpha_{n},\beta_{1},\beta_{2},\ldots, \beta_{n}}(x)\left[\psi(x+2)-\sum_{i=1}^n\alpha_i \gamma_i \psi(\alpha_i x+\beta_i)\right] 
\end{equation*}
Under the given conditions that $\min (\alpha_i, \gamma_i) \geq 1$ and $\beta_i \geq 2$ and as digamma ($\psi$) is an increasing function, we deduce that $\sum_{i=1}^n \alpha_i \gamma_i\psi(\alpha_i x+\beta_i) \geq \psi(x+2)$ and consequently the function $f^{\gamma_{1},\gamma_{2},\ldots,\gamma_{n}}_{\alpha_{1},\alpha_{2},\ldots,\alpha_{n},\beta_{1},\beta_{2},\ldots, \beta_{n}}(x)$ is decreasing on $[1, \infty)$. This implies that the sequence $\left(y_{k}\right)_{k \geq 1}$ monotonically decreases for all $\alpha_i \geq 1, \beta_i \geq 2$ and $\gamma_i \geq 1$ for all $i =1,2,\cdots,n$. Therefore

\begin{equation*}
\left|\left(\mathbf{F}(z)\right)^{\prime}-1\right|<\sum_{k=1}^{\infty} \frac{y_{1}}{k!}=y_{1}(e-1) 
\end{equation*}
This implies that

$$
\left|\left(\mathbf{F}(z)\right)^{\prime}-1\right|<1, \quad z \in \mathcal{D}
$$
Hence, the function $\mathbf{F}(z)$ is convex in $\mathcal{D}_{1 / 2}$ by Lemma \ref{convexnew}. This completes the proof.
    \end{enumerate}
\end{proof}

\begin{example}
Let $n=3$ and consider the parameter sets
$(\alpha_1,\beta_1,\gamma_1)=(1.2,\,1.0,\,1.0)$,
$(\alpha_2,\beta_2,\gamma_2)=(1.5,\,1.0,\,1.2)$,
$(\alpha_3,\beta_3,\gamma_3)=(1.0,\,1.5,\,1.2)$.
Here $\min(\alpha_i,\gamma_i)\ge1$, $\beta_i>0$, $\alpha_2+\beta_2=2.5\ge2$, and
$\prod_{i=1}^{3}\Gamma(\beta_i)^{\gamma_i}
 =\Gamma(1.5)^{1.2}\approx0.86$.
Since
$\prod_{i=1}^{3}[\Gamma(\alpha_i+\beta_i)]^{\gamma_i}
 =\Gamma(2.2)\cdot\Gamma(2.5)^{1.2}\cdot\Gamma(2.5)^{1.2}\approx2.18$,
the condition $(e-1)\prod\Gamma(\beta_i)^{\gamma_i}<\prod[\Gamma(\alpha_i+\beta_i)]^{\gamma_i}$
is satisfied, so by part~(a) of the preceding theorem, $\mathbf{F}$ is starlike in $\mathcal{D}_{1/2}$. In Figure~\ref{dhalf_starlike}, the orange region is the image of $\mathcal{D}_{1/2}$ under $\mathbf{F}(z)$. The image encircles the origin as a simple closed curve, confirming that $\mathbf{F}$ satisfies the starlikeness condition on $\mathcal{D}_{1/2}$.
\begin{figure}[h!]
    \centering
    \includegraphics[width=0.5\linewidth]{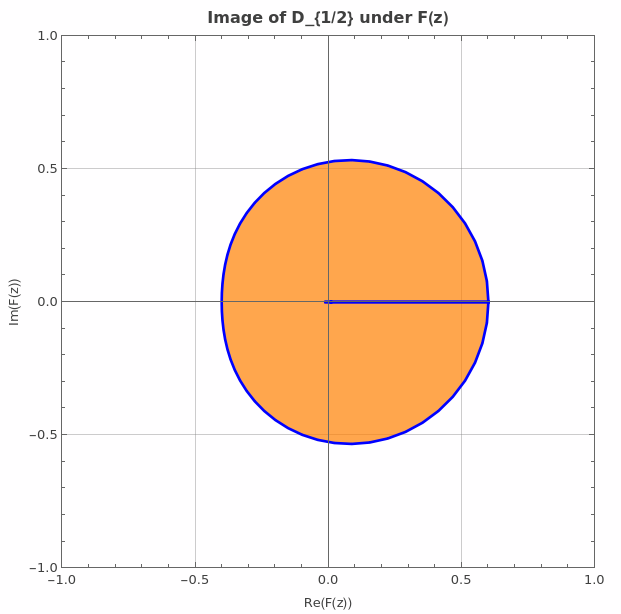}
    \caption{$n=3$ with $(\alpha_1,\beta_1,\gamma_1)=(1.2,1.0,1.0)$,
    $(\alpha_2,\beta_2,\gamma_2)=(1.5,1.0,1.2)$,
    $(\alpha_3,\beta_3,\gamma_3)=(1.0,1.5,1.2)$: image of $\mathcal{D}_{1/2}$ under $\mathbf{F}(z)$ (orange) with boundary (blue)}
    \label{dhalf_starlike}
\end{figure}
\end{example}

\begin{example}
Let $n=3$ and consider the parameter sets
$(\alpha_1,\beta_1,\gamma_1)=(1.0,\,2.0,\,1.2)$,
$(\alpha_2,\beta_2,\gamma_2)=(1.5,\,2.0,\,1.0)$,
$(\alpha_3,\beta_3,\gamma_3)=(1.2,\,2.5,\,1.0)$.
Here $\min(\alpha_i,\gamma_i)\ge1$, $\beta_i\ge2$ for all $i$, $\alpha_2+\beta_2=3.5\ge2$, and
$\prod_{i=1}^{3}\Gamma(\beta_i)^{\gamma_i}
 =\Gamma(2.0)^{1.2}\cdot\Gamma(2.0)\cdot\Gamma(2.5)\approx1.33$.
Since
$\prod_{i=1}^{3}[\Gamma(\alpha_i+\beta_i)]^{\gamma_i}
 =\Gamma(3.0)^{1.2}\cdot\Gamma(3.5)\cdot\Gamma(3.7)\approx31.84$,
the condition $2(e-1)\prod\Gamma(\beta_i)^{\gamma_i}<\prod[\Gamma(\alpha_i+\beta_i)]^{\gamma_i}$
is satisfied, so by part~(b) of the preceding theorem, $\mathbf{F}$ is convex in $\mathcal{D}_{1/2}$. In Figure~\ref{dhalf_convex}, the orange region is the image of $\mathcal{D}_{1/2}$ under $\mathbf{F}(z)$. 
\begin{figure}[h!]
    \centering
    \includegraphics[width=0.5\linewidth]{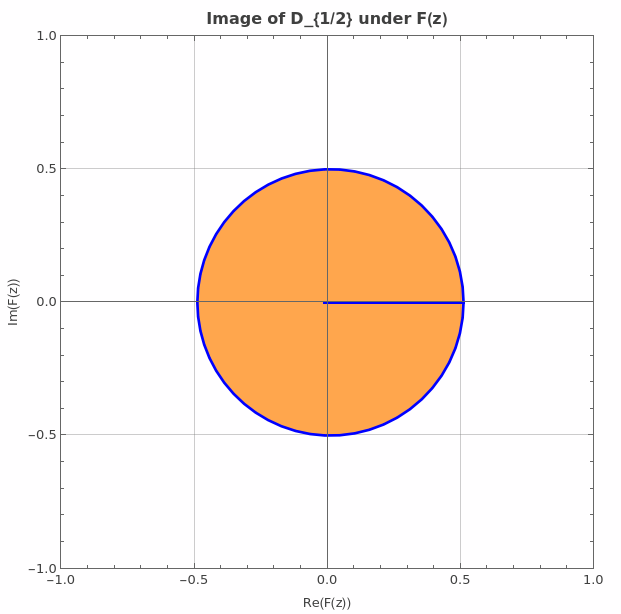}
    \caption{$n=3$ with $(\alpha_1,\beta_1,\gamma_1)=(1.0,2.0,1.2)$,
    $(\alpha_2,\beta_2,\gamma_2)=(1.5,2.0,1.0)$,
    $(\alpha_3,\beta_3,\gamma_3)=(1.2,2.5,1.0)$: image of $\mathcal{D}_{1/2}$ under $\mathbf{F}(z)$ (orange) with boundary (blue)}
    \label{dhalf_convex}
\end{figure}
\end{example}
Comparing Figure~\ref{dhalf_convex} with Figure~\ref{dhalf_starlike} we observe that the convex case produces a more rounded image, reflecting the stronger convexity condition.

\section{Conclusion}
In this paper, we have established sufficient conditions for exponential subordination, exponential starlikeness, exponential convexity, starlikeness, convexity, and close-to-convexity of the Le Roy type Mittag-Leffler functions and the multivariate Le Roy type Mittag-Leffler functions on the unit disk $\mathcal{D}$. The results identify parameter regimes involving $\min(\alpha_i,\gamma_i)\ge1$, $\beta_i>0$, and suitable bounds on $\alpha_j+\beta_j$ that ensure these geometric properties. The close-to-convexity result is notable as it requires only the minimal hypothesis $\alpha_j+\beta_j\ge2$ together with the natural restrictions $\min(\alpha_i,\gamma_i)\ge1$ and $\beta_i>0$ and no additional coefficient inequality.

 The results were proved using a combination of ratio estimates, digamma bounds, and classical subordination criteria, and were illustrated graphically through images of the unit disk generated in Wolfram Mathematica. Our work also generalizes the known geometric properties of three-parameter Le Roy type Mittag-Leffler functions to the multivariate $3n$-parameter setting, providing a unified framework for studying the geometric characteristics of this class of special functions. These results contribute to the geometric function theory of Le Roy type functions, opening directions for further investigation such as extending the framework to lemniscate starlikeness and convexity, or exploring analogous properties for other multi-index special functions.

\end{document}